\input amstex 
\font\we=cmb10 at 14.4truept
\font\li=cmb10 at 12truept
\documentstyle{amsppt} 
\catcode`\@=11
\redefine\logo@{}
\catcode`\@=13
\hsize=6.5truein
\vsize=9.5truein
\topmatter 
\title\nofrills {\we Non-Abelian  L-Functions For Number Fields} 
\endtitle 
\author {\li Lin WENG} 
\endauthor 
\affil {\bf Graduate School of Mathematics, Kyushu University, Japan}\endaffil
%\address 6-10-1 Hakozaki, Fukuoka, Japan\endaddress
\NoRunningHeads 
\endtopmatter 
\TagsOnRight
In this paper we introduce non-abelian zeta functions and more generally
non-abelian $L$-functions for number fields, based on 
geo-arithmetical cohomology, geo-arithmetical truncation and Langlands' theory
of Eisenstein series. More precisely, in Chapter I, we start with a new yet 
natural geo-arithmetical cohomology and a geo-arithmetical stability in order 
to define genuine non-abelian zeta functions for number fields. Then, 
in Chapter II, using examples for the field of rationals, we explain a
relation between these non-abelian zeta functions and Eisenstein series,
and point out where non-abelian contributions of such non-abelian zetas
come. After that, in Chapter III, first, we give a more general 
geo-arithmetical truncation, and compare it with Arthur's analytic truncation. 
Then,  we  
define our general non-abelian $L$-functions as the integration of Eisesntein
series associated to $L^2$-automorphic forms over the moduli spaces obtained 
from just mentioned general geo-arithmetic  truncations, and to expose basic
properties for these non-abelian $L$-functions. In Chapter IV, we, 
along with the line of Rankin-Selberg(-Langlands-Arthur) method, calculate 
the so-called abelian part of our non-abelian 
$L$-functions. Finally, in an Appendix, we establish a non-abelian class 
field theory for function fields over complex numbers to justify our approach 
to non-abelian Reciprocity Law.
\vskip 0.45cm
\centerline {\li Chapter I. Non-Abelian Zeta Functions}
\vskip 0.30cm
\centerline {\li I.1. A Geo-Arithmetical Cohomology}
\vskip 0.30cm
\noindent
{\bf (I.1.1) Cohomology Groups}
\vskip 0.30cm
Let $F$ be a number field with discriminant $\Delta_F$. Denote its
(normalized) absolute values by
$S_F$, and write $S_F=S_{\text{fin}}\buildrel \cdot\over\cup S_\infty$,
where $S_\infty$ denotes the collection of all archimedean valuations.
For simplicity, we often use $v$ (resp. $\sigma$) to denote elements in 
$S_{\text{fin}}$ (resp. $S_\infty$).

Denote by ${\Bbb A}={\Bbb A}_F$ the ring of adeles of $F$, by 
$\text{GL}_r({\Bbb A})$ the rank $r$ general linear group over 
${\Bbb  A}$, and write
${\Bbb  A}:={\Bbb  A}_{\text{fin}}\oplus {\Bbb  A}_\infty$ and 
$\text{GL}_r({\Bbb  A}):=\text{GL}_r({\Bbb  A})_{\text{fin}}\times 
\text{GL}_r({\Bbb  A})_\infty$
according to their finite and infinite parts.

For any $g=(g_{\text{fin}}:g_\infty)=(g_v;g_\sigma)\in \text{GL}_r({\Bbb  A})$,
define the injective morphism $i(g):=i(g_\infty):F^r\to 
{\Bbb  A}^r$ by $f\mapsto (f;g_\sigma\cdot f)$. Let $F^r(g):=
\text{Im}\big(i(g)\big)$ and set
$$\aligned {\Bbb  A}^r(g):=\endaligned
\left\{\aligned (a_v;a_\sigma)\in {\Bbb A}^r:\endaligned 
\aligned \text{(i)}&\ g_v(a_v)\in {\Cal  O}_v^r,\ \forall v;\\ 
\text{(ii)}&\ \exists f\in F^r\ \text{such that}\\
\text{(ii.a)}&\ g_v(f)\in {\Cal 
O}_v^r,\,\forall v;\quad \text{(ii.b)}\ (f;a_\sigma)=i(g_\infty)(f) \endaligned\right\}$$ Then we
have the following 9-diagram with exact columns and rows:
$$\matrix &&0&&0&&0&&\\
&&\downarrow&&\downarrow&&\downarrow&&\\
0&\to&{\Bbb A}^r(g)\cap F^r(g)&\to&{\Bbb A}^r(g)&\to&{\Bbb A}^r(g)/
{\Bbb A}^r(g)\cap F^r(g)&\to&0\\
&&\downarrow&&\downarrow&&\downarrow&&\\
0&\to& F^r(g)&\to&{\Bbb A}^r&\to&{\Bbb A}^r/F^r(g)&\to&0\\
&&\downarrow&&\downarrow&&\downarrow&&\\
0&\to& F^r(g)/{\Bbb A}^r(g)\cap F^r(g)&\to&{\Bbb A}^r/{\Bbb A}^r(g)&\to&
{\Bbb A}^r/{\Bbb A}^r(g)+F^r(g)&\to&0\\
&&\downarrow&&\downarrow&&\downarrow&&\\
&&0&&0&&0.&&\endmatrix$$

Motivated by this and Weil's adelic cohomology theory for divisors over
algebraic curves, (see e.g, [Se2] and [W2]), we introduce the
following
\vskip 0.30cm
\noindent
{\bf Definition.} {\it For any $g\in \text{GL}_r({\Bbb A})$, define its
0-th and 1-st geo-arithmetical cohomology groups by
$$H^0({\Bbb A}_F,g):={\Bbb A}^r(g)\cap F^r(g),\qquad\text{and}\qquad 
H^1({\Bbb A}_F,g):={\Bbb A}^r/{\Bbb A}^r(g)+F^r(g).$$}

\noindent
{\bf Proposition 1.} (Serre Duality=Pontrjagin Duality) {\it As locally compact
groups, $$H^1({\Bbb A}_F,g)\simeq 
\widehat{H^0({\Bbb A}_F,\kappa_F\otimes g^{-1})}.$$ Here $\kappa_F$
denotes an idelic dualizing element of $F$, and $\hat\cdot$  the
Pontrjagin dual. In particular,
$H^0$ is  discrete and $H^1$ is compact.}
\vskip 0.30cm
\noindent
{\it Remark.} For $v\in S_{\text{fin}}$, denote by $\partial_v$ the local
different of $F_v$, the $v$-completion of $F$, and by
${\Cal  O}_v$ the valuation ring with $\pi_v$ a local parameter. Then
$\partial_v=:\pi_v^{\text{ord}_v(\partial_v)}\cdot {\Cal  O}_v$. We call
$\kappa_F:=(\partial_v^{\text{ord}_v(\partial_v)};1)\in {\Bbb I}_F:=
\text{GL}_1({\Bbb A})$ 
an idelic dualizing element of $F$.
\vskip 0.30cm
\noindent
{\it Proof.} As usual, fix a basic character $\chi$ on ${\Bbb A}$ by 
$\chi=:(\chi^{(r)}_v;\chi^{(r)}_\sigma)$ where $\chi_v:=\lambda_v\circ 
\text{Tr}^{F_v}_{{\Bbb Q}_v}$ with $\lambda_v:{\Bbb Q}_v\to {\Bbb Q}_v/
{\Bbb Z}_v\hookrightarrow {\Bbb Q}/{\Bbb Z}\hookrightarrow {\Bbb R}/{\Bbb Z}$, 
and $\chi_\sigma:=\lambda_\infty\circ \text{Tr}_{\Bbb R}^{F_\sigma}$ with
$\lambda_\infty:{\Bbb R}\to {\Bbb R}/{\Bbb Z}$.  Then the pairing
$(x,y)\mapsto e^{2\pi i\chi(x\cdot y)}$ induces natural isomorphisms 
$\widehat {{\Bbb A}^r}\simeq {\Bbb A}^r$ (as
locally compact groups) and $(F^r)^\perp\simeq F^r$ (as discrete subgroups). 
With this, a direct local calculation, see e.g, [Ta], shows that
 $\big({\Bbb A}^r(g)\big)^\perp={\Bbb A}^r(\kappa_F\otimes g^{-1})$ and
$\big(F^r(g)\big)^\perp\simeq F^r(\kappa_F\otimes g^{-1})$. This 
completes the proof since
$$\Big({\Bbb A}^r(g)\cap F^r(g)\Big)^\perp=\big({\Bbb A}^r(g)\big)^\perp
+\big(F^r(g)\big)^\perp.$$
\vskip 0.30cm
\noindent
{\bf (I.1.2) Geo-Arithmetical Counts}
\vskip 0.30cm
Motivated by the Pontrjagin duality and the fact that the dimension of a vector space is
equal to the dimension of its dual, one  basic principal we
adopt in counting locally compact groups is the following:
\vskip 0.30cm
\noindent
{\bf Counting Axiom.} {\it If $\#_{\text{ga}}$  counts a certain
class of locally compact groups $G$,  then $\#_{\text{ga}}(G)=
\#_{\text{ga}}(\hat G)$.}
\vskip 0.30cm
Practically, our counts of arithmetic cohomology groups are based on the 
Fourier inverse
formula, or more accurately, the Plancherel formula in Fourier analysis 
over locally compact
groups. (See e.g, [Fo].) 

While  any reasonable test function on ${\Bbb A}^r$ would do,  as a
continuation of a more traditional mathematics and also for simplicity, 
we set $f:=\prod_v
f_v\cdot\prod_\sigma f_\sigma$, where 

\noindent
(i) $f_v$ is the characteristic function of ${\Cal  O}_v^r$;

\noindent
(ii.a) $f_\sigma(x_\sigma):=e^{-\pi|x_\sigma|^2/2}$ if $\sigma$ is real; and

\noindent
(ii.b) $f_\sigma(x_\sigma):=e^{-\pi|x_\sigma|^2}$ if $\sigma$ is complex.

\noindent
Moreover, we take the following  normalization for the Haar measure
$dx$, which we call {\it standard}, on ${\Bbb A}$: locally 

for $v$, $dx$ is the measure for which
${\Cal  O}_v$ gets the measure $N(\partial_v)^{-1/2}$;

while for $\sigma$ real (resp. complex), $dx$ is the ordinary Lebesgue measure
(resp. twice the ordinary Lebesgue measure).
\vskip 0.30cm
\noindent
{\bf Definition.} {\it (1) The geo-arithmetical counts of the 0-th and the
1-st  cohomology groups for  $g\in \text{GL}_r({\Bbb A})$ are 
defined to be
$$\eqalign{\#_{\text{ga}}\big(H^0({\Bbb A}_F,g)\big):=
&\#_{\text{ga}}\Big(H^0({\Bbb A}_F,g); f,
dx\Big):=\int_{H^0({\Bbb A}_F,g)}|f(x)|^2dx;\cr
\#_{\text{ga}}\big(H^1({\Bbb A}_F,g)\big):=
&\#_{\text{ga}}\Big(H^1({\Bbb A}_F,g); \hat f,
d\xi\Big):=\int_{H^1({\Bbb A}_F,g)}|\hat f(\xi)|^2d\xi.\cr} $$ 
Here $dx$ denotes 
(the restriction of) the standard Haar measure on ${\Bbb A}$,  $d\xi$
(the induced quotient measure from) the dual measure (with
respect to $\chi$), and
$\hat f$  the corresponding Fourier transform of $f$;
\vskip 0.30cm
\noindent
(2) The 0-th and the 1-st geo-arithmetical cohomologies of 
$g\in \text{GL}_r({\Bbb A})$ are defined to be
$$h^0({\Bbb A}_F,g):=\log \Big(\#_{\text{ga}}
\big(H^0({\Bbb A}_F,g)\big)\Big)\quad
\text{and}\quad h^1({\Bbb A}_F,g):=
\log \Big(\#_{\text{ga}}\big(H^1({\Bbb A}_F,g)\big)\Big).$$}
\vskip 0.30cm
\noindent
{\it Remark.} While the groups $H^0({\Bbb A}_F,g)$ and $H^1({\Bbb A}_F,g)$
do depend on $g\in \text{GL}_r({\Bbb A})$, the counts $h^0({\Bbb A}_F,g)$ and 
$h^1({\Bbb A}_F,g)$ depend only on the class of
$g\in \text{GL}_r(F)\backslash \text{GL}_r({\Bbb A})$ by the product formula.
\vskip 0.30cm
\noindent
{\bf (I.1.3) Serre Duality and Riemann-Roch}
\vskip 0.30cm
For the arithmetic cohomologies just introduced, we have the following
\vskip 0.30cm
\noindent
{\bf Theorem.} For any $g\in \text{GL}_r(F)\backslash \text{GL}_r({\Bbb A})$,

\noindent
(1) (Serre Duality) $h^1({\Bbb A}_F,g)=
h^0({\Bbb A}_F,\kappa_F\otimes g^{-1})$;

\noindent
(2) (Riemann-Roch Theorem) 
$$h^0({\Bbb A}_F,g)-h^1({\Bbb A}_F,g)=\text{deg}(g)-{r\over 2}\cdot\log|\Delta_F|.$$

\noindent
{\it Proof.} By the choice of our $f$ in the definition, (1) is a direct 
consequence of the
topological  Serre duality, i.e, Prop. 1.1, and the Plancherel Formula;
 while (2) is a direct consequence of the Serre duality just proved and Tate's
Riemann-Roch theorem ([Ta, Thm. 4.2.1] and/or [Neu, XIV, \S6]), i.e, the Poisson summation
formula, by the fact that
$\Big(H^0({\Bbb A}_F,g)\Big)^\perp =H^0({\Bbb A}_F,\kappa_F\otimes g^{-1})$.
This then completes the proof.
\vskip 0.30cm
Often, for our own convenience, we also write 
$e^{h^i({\Bbb A}_F,g)}:=H^i_{\text{ga}}(F,g),\ i=1,2$. With this, the above 
additive version may be rewritten as 
\vskip 0.30cm
\noindent
{\bf Theorem}$'$. For any $g\in \text{GL}_r(F)\backslash 
\text{GL}_r({\Bbb A})$,

\noindent
(1) (Serre Duality) $H^1_{\text{ga}}(F,g)=H^0_{\text{ga}}
(F,\kappa_F\otimes g^{-1})$;
\vskip 0.30cm
\noindent
(2) (Riemann-Roch Theorem) 
$H^0_{\text{ga}}(F,g)=H^1_{\text{ga}}(F,g)\cdot N(g)\cdot N(\kappa_F)^{-{r\over 2}}$, where as
usual $N(g)$ denotes $e^{{\text{deg}}(g)}$.
\vskip 0.30cm
\noindent
{\it Remarks.} (1) Our work here is motivated by the works of Weil [We2], 
Tate[Ta], Iwasawa [Iw1,2], van der Geer-Schoof [GS], and Li [Li], as well as the works of Lang [L1,2], 
Arakelov [L3],
Szpiro [L3], Parshin [Pa],  Moreno [Mo], Neukirch [Neu],  Deninger [De 1,2], Connes [Co], and  
Borisov [Bor]. For details,
please see the references below, in particular [We2]. 
Also, it would be extremely
interesting if one could relate the work here with that of Connes [Co] and 
Deninger [De1,2,3].

\noindent
(2)  One may apply the discussion here to  wider classes of (multiplicative) 
characters
and test functions so as to obtain a more general class of $L$-functions as
done in Tate's Thesis. But we later would take a substentially different 
approach.

\noindent
(3) Instead of working on $S_F$, the set of equivalence classes of valuations 
of $F$, or better, normalized valuations over $F$, one 
may also work on the set of {\it all} valuations of $F$, which in particular 
is parametrized by $S_F\times {\Bbb R}$ to get a much more refined 
geo-arithmetical cohomology theory. For details, please see [We3].

\vskip 0.8cm
\centerline {\li I.2. New Non-Abelian Zeta Functions}
\vskip 0.30cm
\noindent
{\bf (I.2.1) Intersection Stability}
\vskip 0.30cm
For a metrized vector sheaf $({\Cal  E},\rho)$ on $\text{Spec}({\Cal  O}_F)$, 
define its associated $\mu$-invariant by 
$$\mu({\Cal  E},\rho):={{{\text {deg}}_{\text {Ar}}({\Cal  E},\rho)}\over
{{\text {rank}}({\Cal  E})}},$$ where ${\Cal  O}_F$ denotes the ring of integers 
of a number field $F$ and ${\text {deg}}_{\text {ar}}$  the Arakelov degree of 
$({\Cal  E},\rho)$. (See e.g. [16].) By definition, a proper sub metrized 
vector sheaf $({\Cal  E}_1,\rho_1)$ of $({\Cal  E},\rho)$
consists of a proper sub vector sheaf ${\Cal  E}_1$ of ${\Cal  E}$ such that 
$\rho_1$ is induced from the restriction of $\rho$ via the injection 
${\Cal  E}_1\hookrightarrow {\Cal  E}$.
\vskip 0.30cm
\noindent
{\bf Definition.} {\it A metrized vector sheaf $({\Cal  E},\rho)$ is called 
stable (resp. semi-stable) if for all proper sub metrized vector sheaf 
$({\Cal  E}_1,\rho_1)$ of $({\Cal  E},\rho)$,
$$\mu({\Cal  E}_1,\rho_1)\,<\,\mu({\Cal  E},\rho)\qquad(resp.\qquad 
\mu({\Cal  E}_1,\rho_1)\leq\mu({\Cal  E},\rho)).$$}

\noindent
{\it Remarks.} (1) Despite that we define it independently,
the intersection stability in arithmetic, motivated by Mumford's work [Mu] 
in geometry, was first introduced by Stuhler in [St1,2],  see also 
[Gr1,2], [Mor] and [Bos]. 

\noindent
(2) The intersection stability plays a key role in our work on non-abelian 
class field theory for Riemann surfaces in [We1], appeared as what I call the
{\it micro reciprocity law}. Motivated by this, as a 
fundamental problem, we ask
whether a similar micro reciprocity law holds for number fields. This then
leads naturally to global representations, or better, adelic representations,
of Galois groups. For details about such a discussion, please see [We2].
\vskip 0.30cm
Motivated by Weil's adelic interpretation of vector bundles over 
function fields, for $g=(g_v;g_\sigma)\in \text{GL}_r({\Bbb A})$, introduce a 
torsion-free ${\Cal  O}_F$-module
$$H^0({\Bbb A}_F,g)_{\text{fin}}:=H^0(\text{Spec}({\Cal  O}_F),g):=
\{f\in F^r:g_vf\in {\Cal  O}_v^r,\forall v\}$$ in
$F^r$. Denote the associated vector sheaf on $\text{Spec}({\Cal  O}_F)$ by 
${\Cal  E}(g)$, that is,
$${\Cal  E}(g):=\widetilde {H^0(\text{Spec}({\Cal  O}_F),g)}.$$
Moreover, note that $F^r$, via completion, is densely embedded in 
${\Bbb A}_\infty^r$. Thus to
introduce metrics on ${\Cal  E}(g)$ is equivalent to assigning metrics on the 
associated determinants, i.e., on the top exterior products, of the initial 
data. (See e.g, [L3, Chap. V].) Hence without loss of
generality, we may assume that $r=1$. In this case, the metric on 
${\Cal  E}(g)$ (associated to $g$)
is defined to be the one such that for the rational section $1\in F$,
$$\|1\|_\sigma:=\|g_\sigma\|_\sigma:=
|g_\sigma|^{N_\sigma:=[F_\sigma:{\Bbb Q}_\sigma]}.$$ Denote
such a metric on ${\Cal  E}(g)$ by $\rho(g)$ for $g\in \text{GL}_r({\Bbb A})$. 

As such, we obtain a canonical map 
$({\Cal  E}(\cdot),\rho(\cdot)):\text{GL}_r({\Bbb A})\to
\Omega_{\text{Spec}({\Cal  O}_F),r}$ by assigning $g$ to 
$({\Cal  E}(g),\rho(g))$, where $\Omega_{\text{Spec}({\Cal  O}_F),r}$ 
denotes the collection of all isometrized classes of metrized vector sheaves of
rank $r$ over $\text{Spec}({\Cal  O}_F)$. Clearly, 
$({\Cal  E}(\cdot),\rho(\cdot))$ factors through
the quotient group $\text{GL}_r(F)\backslash \text{GL}_r({\Bbb A})$ where 
$\text{GL}_r(F)$ is embedded diagonally in $\text{GL}_r({\Bbb A})$. 
Denote this resulting map  by 
$({\Cal  E}(\cdot),\rho(\cdot))$ too by an abuse of notation. (See Ch.II below.)

Denote by ${\Cal  M}_{F,r}(d)$ the subset of 
$\Omega_{\text{Spec}({\Cal  O}_F),r}$ consisting of the classes resulting from 
semi-stable metrized vector sheaves of (Arakelov) degree $d$. Since for a 
fixed degree, the semi-stability condition is a bounded and closed one, 
with respect to the natural topology, ${\Cal  M}_{F,r}(d)$ is compact. 
(See e.g. below or [Bo1,2], [Gr1] and/or [St1,2].)

Denote by ${\Cal  M}_{{\Bbb A}_F,r}(d)\subset 
\text{GL}_r(F)\backslash \text{GL}_r({\Bbb A})$ the
inverse image of ${\Cal  M}_{F,r}(d)$ under the map
$({\Cal  E}(\cdot),\rho(\cdot))$, and denote
the corresponding map induced from the restriction by
$$\Pi_{F,r}(d):{\Cal  M}_{{\Bbb A}_F,r}(d)\to {\Cal  M}_{F,r}(d)$$ which we call 
the (algebraic) moment map. As a subquotient of $\text{GL}_r({\Bbb A})$, 
${\Cal  M}_{{\Bbb A}_F,r}(d)$ admits a natural topology, the induced one. 
Also, by a general result due to Borel [Bo1], which in
our case is more or less obvious, the fibers of $\Pi_{F,r}(d)$ are all 
compact.  Thus in particular, ${\Cal  M}_{{\Bbb A}_F,r}(d)$, which we call 
the moduli space of semi-stable
adelic bundles of rank $r$ and degree $d$, is compact. Moreover, as a 
subquotient of 
$\text{GL}_r({\Bbb A})$, ${\Cal  M}_{{\Bbb A}_{F,r}}(d)$ carries a natural 
measure induced from the standard one on $\text{GL}_r({\Bbb A})$, 
which we call the {\it Tamagawa measure}, and denote  by
$d\mu_{{\Bbb A}_F,r}(d)$. For the same reason, there is also a natural 
measure on 
${\Cal  M}_{F,r}(d)$, which we have reasons to call the {\it hyperbolic measure}, and denote it by $d\mu_{F,r}(d)$. 

Clearly, the total volumes of ${\Cal  M}_{{\Bbb A}_F,r}(d)$ 
(resp. ${\Cal  M}_{F,r}(d)$) with
respect to $d\mu_{{\Bbb A}_F,r}(d)$ (resp. $d\mu_{F,r}(d)$) are very important 
non-commutative invariants for number fields. 
Note that according to what we call the Bombieri-Vaaler trick [BV], i.e.,  
by multiplying $g$ with
$(1;e^{t_\sigma})$ where $t_\sigma:=N_\sigma\cdot t$ with
$t\in {\Bbb R}$, we obtain a natural isomorphism between 
${\Cal  M}_{{\Bbb A}_F,r}(d)$ and 
${\Cal  M}_{{\Bbb A}_F,r}(d-n\cdot t)$ where $n:=[F:{\Bbb Q}]$. 
(Even though it is an open problem that whether
semi-stability is closed under tensor operation [Bos], the case here in which 
one is of rank 1 is rather obvious.) Consequently, the above volumes
are independent of degrees $d$. Denote them by $W_F(r)$ and $w_F(r)$
respectively. Later we will see that $W_F(r)$ and $w_F(r)$ are the residues at 
$s=1$ of our non-abelian zeta functions. 
\vskip 0.30cm
\noindent
{\bf (I.2.2) Functional Equation: A Formal Calculation}
\vskip 0.30cm
Motivated by Weil's treatment of functional equation in [W2], we give a formal discussion to show the key points of functional equation.

For the number field $F$ with discriminant $\Delta_F$, denote by 
${\Cal  M}_{{\Bbb A}_F,r}$ the
moduli space of semi-stable adelic bundles of rank $r$, 
that is, ${\Cal  M}_{{\Bbb A}_F,r}:=
\buildrel\cdot\over\cup_{N\in {\Bbb R}_+} {\Cal  M}_{{\Bbb A}_F,r}[N]$
where ${\Cal  M}_{{\Bbb A}_F,r}[N]:={\Cal  M}_{{\Bbb A}_F,r}(\log\,N)$. 
By using the Bombieri-Vaaler trick in (2.1), as topological spaces, 
${\Cal  M}_{{\Bbb A}_F,r}\simeq 
{\Cal  M}_{{\Bbb A}_F,r}[|\Delta_F|^{r\over 2}]\times {\Bbb R}_+$.
Hence we obtain a natural measure $d\mu$ on ${\Cal  M}_{{\Bbb A}_F,r}$ from the 
Tamagawa 
measures on ${\Cal  M}_{{\Bbb A}_F,r}[N]$ and ${{dT}\over T}$ on ${\Bbb R}_+$. 
 
Recall that for any $E\in {\Cal  M}_{{\Bbb A}_F,r}$, 
 $H_{\text{ga}}^i(F,E):=H_{\text{ga}}^i(F,g)
=e^{h^i({\Bbb A}_F,g)}$ where $g\in \text{GL}_r({\Bbb A})$ is such that 
$E=[g]$. Since for any 
$a\in \text{GL}_r(F)$, $H_{\text{ga}}^i(F,a\cdot g)=H^i_{\text{ga}}(F,g)$ 
by the product formula. $H_{\text{ga}}^i(F,E)$ is well-defined for $i=0,\,1$.
\vskip 0.30cm
With respect to fixed real constants $A,B,C,\alpha$ and $\beta$, introduce the 
formal
integration $Z_{F,r;A,B,C;\alpha,\beta}(s)$ as follows:
$$Z_{F,r;A,B,C;\alpha,\beta}(s)
:=\big(|\Delta_F|^{-{{rB}\over 2}}\big)^s
\int_{E\in {\Cal  M}_{{\Bbb A}_F,r}}\Big(H_{\text{ga}}^0(F,E)^A\cdot 
N(E)^{Bs+C}-N(E)^{\alpha s+\beta}\Big)d\mu(E).$$ 
Then formally,
$$Z_{F,r;A,B,C;\alpha,\beta}(s)=I(s)-II(s)+III(s),$$ where 
$$\eqalign{I(s):=&\big(|\Delta_F|^{-{{rB}\over 2}}\big)^s
\int_{E\in {\Cal  M}_{{\Bbb A}_F,r},N(E)\leq |\Delta_F|^{r\over 2}}
\Big(H_{\text{ga}}^0(F,E)^A\cdot
N(E)^{Bs+C}-N(E)^{\alpha s+\beta}\Big)d\mu(E);\cr
II(s):=&\big(|\Delta_F|^{-{{rB}\over 2}}\big)^s
\int_{E\in {\Cal  M}_{{\Bbb A}_F,r},N(E)\geq |\Delta_F|^{r\over 2}}
N(E)^{\alpha s+\beta}d\mu(E);\cr
III(s):=&\big(|\Delta_F|^{-{{rB}\over 2}}\big)^s
\int_{E\in {\Cal  M}_{{\Bbb A}_F,r},N(E)\geq |\Delta_F|^{r\over 2}}
H_{\text{ga}}^0(F,E)^A\cdot N(E)^{Bs+C} d\mu(E).\cr} $$
By Theorem$'$ 1.3, i.e., the multiplicative version of 
Serre duality and Riemann-Roch theorem, we have
$$III(s)=
\big(|\Delta_F|^{-{{rB}\over 2}}\big)^{-s-{{A+2C}\over B}}
\int_{E\in {\Cal  M}_{{\Bbb A}_F,r},N(E)\leq |\Delta_F|^{r\over 2}}
H_{\text{ga}}^0(F,E)^A\cdot N(E)^{B(-s-{{A+2C}\over B})+C} d\mu(E),$$ 
by the fact that $N(E_1\otimes E_2^\vee)=
N(E_1)^{\text{rank}(E_2)}\cdot N(E_2)^{-\text{rank}(E_1)}.$ Hence, formally,
$$Z_{F,r;A,B,C;\alpha,\beta}(s)=I(s)+I(-s-{{A+2C}\over B})-II(s)+IV(s),$$ where
$$IV(s):=\big(|\Delta_F|^{-{{rB}\over 2}}\big)^{-s-{{A+2C}\over B}}
\int_{E\in {\Cal  M}_{{\Bbb A}_F,r},N(E)\leq |\Delta_F|^{r\over 2}}
N(E)^{\alpha\cdot(-s-{{A+2C}\over B})+\beta}d\mu(E).$$

Moreover, by definition,
$$\eqalign{-II(s)
=&-\int_{E\in {\Cal  M}_{{\Bbb A}_F,r},
N(E)\geq \Delta_F^{r\over 2}} 
N(E)^{\alpha s+\beta}d\mu(E)
=-\int_{{\Cal  M}_{{\Bbb A}_F,r}[|\Delta_F|^{r\over 2}]} d\mu(E)
\cdot\int_1^\infty t^{\alpha t+\beta}{{dt}\over t}\cr
=&-W_F(r)\cdot {{t^{\alpha s+\beta}}\over
{\alpha s+\beta}}\Big|_1^\infty=W_F(r)\cdot {1\over {\alpha s+\beta}},\cr} $$ 
provided that $\alpha s+\beta<0$.

Similarly,
$$IV(s)=\int_{{\Cal  M}_{{\Bbb A}_F,r}[|\Delta_F|^{r\over 2}]}d\mu(E)
\cdot\int_0^1 t^{\alpha(-s-{{A+2C}\over B})+\beta}{{dt}\over t}
=W_F(r)\cdot {1\over {\alpha(-s-{{A+2C}\over B})
+\beta}}$$ provided that $\alpha(-s-{{A+2C}\over B}) +\beta>0$.

Therefore, formally, 
$$Z_{F,r;A,B,C;\alpha,\beta}(s)=I(s)+I\big(-s-{{A+2C}\over B}\big)
+W_F(r)\cdot\Big({1\over{\alpha s+\beta}}+{1\over {\alpha(-s-{{A+2C}\over
B})+\beta}}\Big).$$
As a direct consequence, we have the following
\vskip 0.30cm
\noindent
{\bf Functional Equation.} {\it With the same notation as above, formally,
$$Z_{F,r;A,B,C;\alpha,\beta}(s)=Z_{F,r;A,B,C;\alpha,\beta}(-s-{{A+2C}\over
B}).$$}
\vskip 0.30cm
\noindent
{\bf (2.3) Non-Abelian Zeta Functions for Number Fields}
\vskip 0.45cm
To justify the arguments in (I.2.2), we consider two types of convergence 
problems.
\vskip 0.30cm
\noindent
{\it Type 1}. Convergence for $II(s)$ and $IV(s)$, where
$$\align II(s)=&\big(|\Delta_F|^{-{{rB}\over 2}}\big)^s
\int_{E\in {\Cal  M}_{{\Bbb A}_F,r},N(E)\geq |\Delta_F|^{r\over 2}}
N(E)^{\alpha s+\beta}d\mu(E);\\
IV(s)=&\big(|\Delta_F|^{-{{rB}\over 2}}\big)^{-s-{{A+2C}\over B}}
\int_{E\in {\Cal  M}_{{\Bbb A}_F,r},N(E)\leq |\Delta_F|^{r\over 2}}
N(E)^{\alpha\cdot(-s-{{A+2C}\over B})+\beta}d\mu(E).\endalign$$
From the calculation in (I.2.2),  when
$\text{Re}(\alpha\cdot s+\beta)<0$ and 
$\text{Re}\big(\alpha\cdot(-s-{{A+2C}\over
B}\big)+\beta\big)>0$, being holomorphic functions,
$$II(s)=-W_F(r)\cdot {1\over {\alpha\cdot s+\beta}}\qquad\text{and}\qquad 
IV(s)=W_E(r)\cdot
{1\over {\alpha\cdot(-s-{{A+2C}\over B}\big)+\beta}}.$$

\noindent
{\it Type 2}. Convergence for $I(s)$ and $I(-s-{{A+2C}\over B})$ where
$$I(s)=
\big(|\Delta_F|^{-{{rB}\over 2}}\big)^s\int_{E\in {\Cal 
M}_{{\Bbb A}_F,r},N(E)\leq |\Delta_F|^{r\over 2}}
\Big(H_{\text{ga}}^0(F,E)^A\cdot
N(E)^{Bs+C}-N(E)^{\alpha s+\beta}\Big)d\mu(E).$$

By the discussion above about $II(s)$, unless $B=\alpha$, $I(s)$ and 
$I(-s-{{A+2C}\over B})$, and
hence $Z_{F,r;A,B,C;\alpha,\beta}(s)$ cannot be meromorphically extended as 
a meromorphic
function to the whole $s$-plane. (See also the discussion below.) Thus, we 
introduce the following
\vskip 0.30cm
\noindent
{\bf Compatibility Conditions}: $\alpha=B$ and $\beta=C$. 
\vskip 0.30cm
With this,  by a change of variables, we may afterwards assume that
$B=1$ and $C=0$ without loss of generality so as to obtain the following 
integration:
$$\align Z_{F,r;A}(s):=&Z_{F,r;A,-1,0;-1,0}(s)\\
:=&
\big(|\Delta_F|^{{{r}\over 2}}\big)^s
\cdot\int_{E\in {\Cal  M}_{{\Bbb A}_F,r}}
\Big(H_{\text{ga}}^0(F,E)^A-1\Big)\cdot N(E)^{-s}d\mu(E).\endalign$$

\noindent
Furthermore, for any $g\in \text{GL}_r({\Bbb A})$, note that 
$H^0({\Bbb A}_F,g)$ is discrete, so from the definition, 
the associated measure is simply the
standard counting measure. Thus by writing down each term precisely, 
$$\align H_{\text{ga}}^0(F,g)=&1+\sum_{\alpha\in H^0(\text{Spec}{\Cal  O}_F,
{\Cal  E}(g))\backslash
\{0\}}\exp\Big(-\pi\sum_{\sigma:{\Bbb R}}|g_\sigma\cdot\alpha|^2-
2\pi\sum_{\sigma:{\Bbb C}}|g_\sigma\cdot \alpha|^2\Big)\\
=:&
1+H_{\text{ga}}'(F,g).\endalign$$ In this expression,
the first term is simply the constant function 1 on the moduli space, 
while each term in the second
decays exponentially. Hence, by a standard argument about convergence of an 
integration of theta
series for higher rank lattices, see. e.g, [Neu], and the
fact that 1 in the first term $H_{\text{ga}}^0(F,E)^A$ cancels with 1
appeared as the 
second term in the
combination $H_{\text{ga}}^0(F,E)^A-1$, we conclude that
$II(s)$ and $II(-s-A)$ are all holomorphic functions, provided $A>0$. 
All in all,
we have proved the following
\eject
\vskip 0.30cm
\noindent
{\bf Main Theorem A.} {\it For any strictly positive real number $A$, 
$$Z_{F,r;A}(s)=:\big(|\Delta_F|^{{{r}\over 2}}\big)^s\cdot
\int_{E\in {\Cal  M}_{{\Bbb A}_F,r}}\Big(H_{\text ga}^0(F,E)^A-1\Big)
\cdot N(E)^{-s}d\mu(E)$$  
is holomorphic when ${\text Re}(s)>A$. Moreover,

\noindent
(1) $Z_{F,r;A}(s)$ admits a meromorphic continuation to the whole complex 
$s$ plane which has 
only two singularities, simple poles at $s=0$ and $s=A$ with the same 
residue $W_F(r)$, i.e, the Tamagawa
volume of the moduli space ${\Cal  M}_{{\Bbb A}_F,r}[|\Delta_F|^{r\over 2}]$;

\noindent
(2) $Z_{F,r;A}(s)$ satisfies
the functional equation $$Z_{F,r;A}(s)=Z_{F,r;A}(A-s).$$}
\noindent
{\it Remark.} As to be seen in Chapter II, the non-abelian invariant $W_F(r)$
may be written in terms of (combinations of products of) special values of 
abelian zeta functions. One may then view this as one of the main application 
of our non-abelian zeta (and $L$-)functions to the classical abelian 
$L$-functions: We obtain a universal relation among special values of abelian 
$L$-functions with the help of the functional equation 
since $W_F(r)$ may be evaluated geometrically -- 
by the well-known formula of Siegel,
we only need to calculate volumes of the cut-off neighbourhood of cusps, a
relatively easy part if a concrete description of the moduli space is known. 
\vskip 0.30cm
\noindent
{\bf Main Definition A.} {\it The function 
$$\xi_{F,r}(s):=\hat\zeta_{F,r}(s):=Z_{F,r}(s):=Z_{F,r;1}(s)$$ is 
called the rank $r$ non-abelian zeta function of $F$.}
\vskip 0.30cm
As well-known to experts, this latest definition may be justified by 
Iwasawa's ICM talk at MIT [Iw1] due to the fact that essentially $\hat\zeta_{F,1}(s)$
is the completed Dedekind zeta function $\hat\zeta_F(s)$ for $F$.
\vskip 0.5cm
\centerline {\li Chapter II. Motivated Examples}
\vskip 0.30cm
\centerline
{\li  II.1. Epstein Zeta Functions and Non-Abelian Zeta Functions}
\vskip 0.30cm
For simplicity, assume that the  number field involved is the field of 
rationals. A lattice $\Lambda$ over ${\Bbb Q}$ is   semi-stable,
by Riemann-Roch, if for any sublattice $\Lambda_1$ of $\Lambda$,
$$\big(\text{Vol}\,\Lambda_1\big)^{\text{rank}\,\Lambda}\geq 
\big(\text{Vol}\,\Lambda\big)^{\text{rank}\,\Lambda_1}.$$
Denote the moduli space of rank $r$ semistable lattices over ${\Bbb Q}$ by 
${\Cal  M}_{{\Bbb Q},r}$, then the lattice version of {\it
rank $r$ non-abelian zeta function} $\xi_{{\Bbb Q},r}(s)$ of ${\Bbb Q}$ is
defined to be  
$$\xi_{{\Bbb Q},r}(s):=
\int_{{\Cal  M}_{{\Bbb Q},r}}\left( e^{h^0({\Bbb Q},\Lambda)}-1\right)
\cdot \big(e^{-s}\big)^{\text{deg}(\Lambda)}\, d\mu(\Lambda),\qquad 
\text{Re}(s)>r,$$
 where $h^0({\Bbb Q},\Lambda):=\log\left(\sum_{x\in \Lambda}
\exp\big(-\pi|x|^2\big)\right)$ and $\text{deg}(\Lambda)$
denotes the Arakelov degree of $\Lambda$. With the discussion in Chap. I,
one checks that

\noindent
(i) $\xi_{{\Bbb Q},1}(s)$ coincides with the (completed)
Riemann-zeta function;

\noindent
(ii) $\xi_{{\Bbb Q},r}(s)$ can be meromorphically extended to the whole 
complex plane;

\noindent
(iii) $\xi_{{\Bbb Q},r}(s)$ satisfies
the functional equation $$\xi_{{\Bbb Q},r}(s)=\xi_{{\Bbb Q},r}(1-s);$$ 

\noindent
(iv) $\xi_{{\Bbb Q},r}(s)$ has only two singularities, simple poles, 
at $s=0,1$, with the same residues
$\text{Vol}\left({\Cal  M}_{{\Bbb Q},r}[1]\right)$, the Tamagawa type volume of 
the space of rank $r$ semi-stable lattice of volume 1. 

Denote by ${\Cal  M}_{{\Bbb Q},r}[T]$ the moduli space of rank $r$ semi-stable 
lattices of volume $T$. We have a 
trivial decomposition $${\Cal  M}_{{\Bbb Q},r}=
\cup_{T>0}{\Cal  M}_{{\Bbb Q},r}[T].$$ Moreover, there is a natural morphism
$${\Cal  M}_{{\Bbb Q},r}[T]\to {\Cal  M}_{{\Bbb Q},r}[1],\qquad \Lambda\mapsto 
T^{1\over r}\cdot\Lambda.$$

With this, $$\align \xi_{{\Bbb Q},r}(s)=&
\int_{\cup_{T>0}{\Cal  M}_{{\Bbb Q},r}[T]}\left( e^{h^0({\Bbb Q},\Lambda)}-1
\right)\cdot \big(e^{-s}\big)^{\text{deg}(\Lambda)}\, d\mu(\Lambda)\\
=&\int_0^\infty T^s{{dT}\over T}
\int_{{\Cal  M}_{{\Bbb Q},r}[1]}\left( 
e^{h^0({\Bbb Q},T^{1\over r}\cdot \Lambda)}-1\right)
\cdot d\mu_1(\Lambda),\endalign$$ where $d\mu_1$ denotes the induced Tamagawa measure 
on ${\Cal  M}_{{\Bbb Q},r}[1]$.

Thus note that $$h^0({\Bbb Q},T^{1\over r}\cdot \Lambda)=
\log\left(\sum_{x\in \Lambda}
\exp\big(-\pi|x|^2\cdot T^{2\over r}\big)\right),$$
and for $B\not=0$,
$$\int_0^\infty e^{-AT^B}T^s{{dT}\over T}={1\over B}\cdot A^{-{s\over B}}
\cdot\Gamma({s\over B}),$$ we have
$$\xi_{{\Bbb Q},r}(s)={r\over 2}\cdot\pi^{-{r\over 2}\, s}
\Gamma({r\over 2}\, s)\cdot\int_{{\Cal  M}_{{\Bbb Q},r}[1]}
\left(\sum_{x\in\Lambda\backslash\{0\}}|x|^{-rs}\right)\cdot d\mu_1(\Lambda).$$
Set now the completed Epstein zeta function, a special kind of Eisenstein 
series, associated to the rank $r$ lattice 
$\Lambda$ over ${\Bbb Q}$ by
$$\hat E(\Lambda;s):=\pi^{-s}\Gamma(s)\cdot \sum_{x\in \Lambda\backslash \{0\}}
|x|^{-2s},$$ then we have the following
\vskip 0.30cm
\noindent
{\bf Proposition.} (Eisenstein series and Non-Abelian 
Zeta Functions) {\it With the same notation as above,
$$\xi_{{\Bbb Q},r}(s)={r\over 2}\int_{{\Cal  M}_{{\Bbb Q},r}[1]}\hat 
E(\Lambda,{r\over 2}s)\,d\mu_1(\Lambda).$$}
 
Thus to study  non-abelian zeta functions, we need to understand Eisenstein 
series and (geo-arithmetical) truncations.
\vskip 0.30cm
\noindent
{\it Remark.} While we may introduce general non-abelian $L$-functions by 
using more general test functions, see e.g, Tate's Thesis where abelian 
version is
discussed. In this paper, we decide to take a different approach using
 Eisenstein series.
(We reminder the reader that for the abelian picture, Eisenstein series 
are not available.)
\vskip 0.30cm
\centerline
{\li II.2. Rankin-Selberg Method: An Example with $SL_2$}
\vskip 0.30cm
From the previous section, we know that 
$$\xi_{{\Bbb Q},2}(s)=\int_{{\Cal  M}_{{\Bbb Q},2}[1]}\hat E(\Lambda,s)\,
d\mu_1(\Lambda).$$
Thus to study $\xi_{{\Bbb Q},2}(s)$, we need to know 

\noindent
(i) what is the moduli space 
of ${\Cal  M}_{{\Bbb Q},2}$? and  

\noindent
(ii) what is the integration of the Eisenstein series $\hat E(\Lambda;s)$ 
over this space.

Before discussing this, let us take a more traditional approach.

Consider the action of $\text{SL}(2,{\Bbb Z})$ on the upper half plane 
${\Cal  H}$.
A standard  \lq fundamental domain\rq$\ $of $SL(2,{\Bbb Z})$ is given by 
$$D=\{z=x+iy\in {\Cal  H}:|x|\leq {1\over 2},y>0,x^2+y^2\geq 1\}.$$ 
Associated to this is aslo the standard completed Eisenstein
series $$\hat E(z;s):=\pi^{-s}\Gamma(s)\cdot\sum_{(m,n)\in {\Bbb Z}^2\backslash \{(0,0)\}}{{y^{s}}\over {|mz+n|^{2s}}}.$$
At this stage, a natural question is to consider the integration
$$\int_D\hat E(z,s){{dx\,dy}\over {y^2}}.\eqno(1)$$
However, as well-known to even graduate students, this integration diverges:
Near the only cusp $y=\infty$, $\hat E(z,s)$ has the  Fourier expansion
$$\hat E(z;s)=\sum_{n=-\infty}^\infty a_n(y,s)e^{2\pi i nx}$$ with 
$$a_n(y,s)=\cases \xi(2s)y^s+\xi(2-2s)y^{1-s},&  \text{if}\ n=0\\
2|n|^{s-{1\over 2}}\sigma_{1-2s}(|n|){\sqrt y}K_{s-{1\over 2}}(2\pi|n|y),& 
\text{if}\ n\not=0, \endcases$$
where $\xi(s)$ denotes the completed Riemann zeta function,
$$\sigma_s(n):=\sum_{d|n}d^s,\quad \text{and}\quad 
K_s(y):={1\over 2}\int_0^\infty e^{-y(t+{1\over t})/2}t^s{{dt}\over t}$$ 
is the K-Bessel function.
Moreover, $$|K_s(y)|\leq e^{-y/2}K_{\text{Re}(s)}(2),\qquad \text{if}\  y>4,
\qquad\text{and}\qquad
K_s=K_{-s},$$ so $a_{n\not=0}(y,s)$ decay exponentially, and the problematic 
term comes from $a_0(y,s)$, which is of slow growth.

Therefore, to make the integration (1) meaningful, we need to cut-off the slow 
growth part. Here we use two ways to do so: 
one is geometric and hence rather direct and simple; while the other is 
analytic, and hence rather technical and traditional, and is dated back to 
Rankin-Selberg.
\vskip 0.30cm
\noindent
(a) {\bf Geometric Truncation}

Draw a horizontal line $y=T\geq 1$ and set
 $$D_T=\{z=x+iy\in D:y\leq T\},\qquad D^T=\{z=x+iy\in D:y\geq T\}.$$ 
Then $D=D_T\cup D^T$.
Introduce the integration
$$I^{\text{Geo}}_T(s):=\int_{D_T}\hat E(z,s)\,{{dx\,dy}\over {y^2}}.$$
Clearly, $I^{\text{Geo}}_T(s)$ is well-defined.
\vskip 0.30cm
\noindent
(b) {\bf Analytic Truncation}

Define a truncated Eisenstein series $\hat E_T(z;s)$ by
$$\hat E_T(z;s):=\cases \hat E(z;s),&\text{if}\ y\leq T\\
\hat E(z,s)-a_0(y;s),&\text{if}\ y>T.\endcases$$
Introduce the integration
$$I_T^{\text{Ana}}(s):=\int_D\hat E_T(z;s)\,{{dx\,dy}\over {y^2}}.\eqno(3)$$
 
With this, from the Rankin-Selberg method, 
one checks that we have the following:
\vskip 0.30cm
\noindent
{\bf Proposition.} (Analytic Truncation=Geometric Truncation in Rank 2) 
{\it With the same notation as above,
$$I_T^{\text{Geo}}(s)=\xi(2s){{T^{s-1}}\over {s-1}}-\xi(2s-1){{T^{-s}}\over {s}}=I_T^{\text{Ana}}(s).\eqno(4)$$}

Each of the above two integrations has its own merit: for the geometric one, 
we keep the Eisenstein series unchanged, while for the analytic one, we keep
the original fundamental domain of ${\Cal  H}$ under $\text{SL}(2,{\Bbb Z})$
as it is.

Note that the nice point about the fundamental domain is that it admits a
modular interpretation. Thus it would be very idealisic if we could at the
same time keep the Eisenstein series unchanged, while offer some integration
domains which appear naturally in certain moduli problems. 
This is exactly what we are seeking for, as indicated by the example 
presented right below. (The general construction is going to be
presented in the last two chapters of this paper.) 
\vskip 0.30cm
\noindent
(c) {\bf  Algebraic Truncation}
\vskip 0.30cm
Now we explain why the above discussion and Rankin-Selberg method have 
anything to do with our non-abelian zeta functions. For this, we introduce 
yet another truncation, the algebraic one.

So back to the moduli space  of rank 2 lattices of volume 1 over ${\Bbb Q}$.
Then classical reduction theory gives a  natural map from this moduli
space to the fundamental domain $D$ above: 
For any lattice $\Lambda$, fix ${\bold x}_1\in \Lambda$ such that its length  
gives the first Minkowski minimum $\lambda_1$ of $\Lambda$.
Then via rotation, we may assume that ${\bold x}_1=(\lambda_1,0)$. 
Further, from the reduction theory 
${1\over {\lambda_1}}\Lambda$ may be viewed as the lattice of the volume
$\lambda_1^{-2}=y_0$ which is generated by $(1,0)$ and $\omega=x_0+iy_0\in D$. 
That is to say, 
the points in $D_T$ are in one-to-one corresponding to
the  rank two lattices of volume one whose
first Minkowski minimum $\lambda_1^{-2}\leq T$, i.e,
$\lambda_1\geq T^{-{1\over 2}}$. Set 
${\Cal  M}_{{\Bbb Q},2}^{\leq {1\over 2}\log T}[1]$  be the moduli space of
rank 2 lattices $\Lambda$ of volume 1 over ${\Bbb Q}$ whose sublattices 
$\Lambda_1$ of rank 1 have degrees 
$\leq {1\over 2}\log T$. As a direct consequence, we have the following
\vskip 0.30cm
\noindent
{\bf Fact.} (Geometric Truncation = Algebraic Truncation) {\it With the same 
notation as above, there is a natural one-to-one, onto morphism 
$${\Cal  M}_{{\Bbb Q},2}^{\leq {1\over 2}\log T}[1]\simeq D_T.$$
In particular, $${\Cal  M}_{{\Bbb Q},2}^{\leq 0}[1]={\Cal  M}_{{\Bbb Q},2}[1]
\simeq D_1.$$}
\vskip 0.30cm
With this, by Proposition I.2, we may then also introduce a much more general type of 
non-abelian zeta functions, parametrized by $T$, with the help of a Harder-Narasimhan
 type discussion on intersection stability. (For details in general, see next chapter.) As a 
special case,  we have the following
\vskip 0.30cm
\noindent
{\bf Main Example.} (Degeneration in Rank 2) {\it With the same notation as 
above, 
$$\xi_{{\Bbb Q},2}(s)=\xi(2s){{1}\over {s-1}}-\xi(2s-1){{1}\over {s}}.
\eqno(5)$$}

Quite disappointed. Isn't it?! After all, what we previously claimed is that
the zeta functions we introduced are non-abelian, yet the calculation
in rank two results only (a combination of) abelian zetas. For me, 
just as happened for many others at the very begining, it took about a half 
year to understand what was involved: As the discussion in the next section
and the last chapter shows all rank two non-abelian zetas degenerate. Despite 
of this, a positive thinking then leads to the following three observations:

\noindent 
(i)  The special values $\zeta(2n)$ and $\zeta(2n-1)$ of the Riemann zeta 
function are intrinsically related by our rank two zeta functions. As 
Eisenstein series may be studied independently,
our formula indicates that non-abelian zetas could be used to understand 
abelian zetas;

\noindent
(ii) The volume of $D_T$ may be evaluated from this formula via a residue 
argument; Moreover, the volume of $D^1$ also have an interpretation in terms 
of the special values of the Riemann zeta.

\noindent
(iii) The dependence on $T$ of the integrations (4) is quite regular: 
The \lq main term' is simply
$$\xi(2s)\cdot {{1}\over {s-1}}-\xi(2s-1)\cdot {1\over {s}}.$$ 
\eject
\vskip 0.30cm
\noindent
\centerline {\li II.3. Where Non-Abelian Contributions Come}
\vskip 0.30cm
In the previous section, we conclude that rank 2 non-abelian zeta function
degenerates. In this section, we explain why this happens and use an example 
of rank 3 zeta functions to point out where exactly the non-abelian 
contributions come. This then justifies why we call the functions introduced 
in the next chapter non-abelian $L$-functions.
\vskip 0.30cm
\noindent
{\bf II.3.1. The Group $SL_3$}
\vskip  0.30cm
The moduli space of rank 3 lattices of volume 1 may be identified with 
$SL(3,{\Bbb Z})\backslash
SL(3,{\Bbb R})/SO(3,{\Bbb R})$.  We start with a description of several 
coordinates for $SL(3,{\Bbb R})/SO(3,{\Bbb R})$.
For this, consider the following standard parabolic subgroups of 
$G=SL(3,{\Bbb R})$.

\noindent
$$P_0=P_{1,1,1}:=\left\{\left(\matrix a_{11}&a_{12}&a_{13}\\
0&a_{22}&a_{23}\\
0&0& a_{33}\endmatrix\right)\in G:a_{ij}\in {\Bbb R}\right\};$$
$$P_1=P_{2,1}:=\left\{\left(\matrix a_{11}&a_{12}&a_{13}\\
a_{21}&a_{22}&a_{23}\\
0&0& a_{33}\endmatrix\right)\in G:a_{ij}\in {\Bbb R}\right\};$$ and
$$P_2=P_{1,2}:=\left\{\left(\matrix a_{11}&a_{12}&a_{13}\\
0&a_{22}&a_{23}\\
0&a_{32}& a_{33}\endmatrix \right)\in G:a_{ij}\in {\Bbb R}\right\}.$$ 

Write the corresponding Langlands decompositions as $P_i=N_iA_iM_i, i=0,1,2$ 
where $N_i$ is the unipotent radical 
of $P_i$, $A_i$ reducible and $M_i$ simple. So, 
$$M_0=\left\{I_3, \left(\matrix 1&0&0\\
0&-1&0\\
0&0& {-1}\endmatrix\right), \left(\matrix -1&0&0\\
0&-1&0\\
0&0&1\endmatrix\right), \left(\matrix
-1&0&0\\
0&1&0\\
0&0&-1\endmatrix\right)\right\}.$$ More generally, if we denote the 
matrices of each subgroup by the corresponding lower-case letters. 
The subgroups above consists of the following elements:
$$\matrix n_0=\left(\matrix 1&a_{12}&a_{13}\\
0&1&a_{23}\\
0&0&1\endmatrix\right);&a_0=\left(\matrix a_{11}&0&0\\
0&a_{22}&0\\
0&0&a_{33}\endmatrix\right);&m_0\in M_0;\\
n_1=\left(\matrix 1&0&x_{1}\\
0&1&t_{1}\\
0&0&1\endmatrix\right);&a_1=\left(\matrix \alpha_{1}&0&0\\
0&\alpha_{1}&0\\
0&0&\alpha^{-1}_{1}\endmatrix\right);&m_1=\left(\matrix *&*&0\\
*&*&0\\
0&0&1\endmatrix\right)\cdot m_0;\\
n_2=\left(\matrix 1&x_2&t_{2}\\
0&1&0\\
0&0&1\endmatrix\right);& a_2=\left(\matrix \alpha^{-2}_{2}&0&0\\
0&\alpha_{2}&0\\
0&0&\alpha_{2}\endmatrix\right);&m_2=\left(\matrix 1&0&0\\
0&*&*\\
0&*&*\endmatrix\right)\cdot m_0,\endmatrix$$
where $a_{ij}, x_i,t_i\in {\Bbb R}, a_{ii}, \alpha_i>0$.

Note that by the Iwasawa decomposition with respect to $P_0$, we 
have $G=A_0^+N_0K$. Thus choose a coset $G/K$ amounts to
choosing an element of $N_0$ and one of $A_0^+$, the identity component 
of $A_0$. Hence, we may identify $G/K$
with $$\left\{Y:=\left(\matrix y_1&0&0\\
0&y_2&0\\
0&0&(y_1y_2)^{-1}\endmatrix \right)\cdot \left(\matrix 1&x_1&x_{2}\\
0&1&x_3\\
0&0&1\endmatrix \right):y_1,y_2>0, x_1,x_2,x_3\in {\Bbb R}\right\}.$$ 
As such it is then convenient to introduce two coordinate systems according 
to the parabolic subgroups $P_1$ and $P_2$, respectively. In fact, 
notice that $M_1/M_1\cap K\simeq SL(2,{\Bbb R})/SO(2,{\Bbb R})$ 
so natural coordinates for $G/K$ are given by
$$\left(\matrix u_1^{1/2}&v_1u_1^{-1/2}&0\\
0&u_1^{-1/2}&0\\
0&0&1\endmatrix \right)\cdot \left(\matrix \alpha_1&0&0\\
0&\alpha_1&0\\
0&0&\alpha_1^{-2}\endmatrix \right)\cdot\left(\matrix 1&0&x_{1}\\
0&1&t_1\\
0&0&1\endmatrix \right)$$ where $z_1=v_1+iu_1$ can be regarded as a point 
in the Poincare upper half plane. Similarly, 
consideration of $P_2$ yields coordinates 
$$\left(\matrix 1&0&0\\
0&u_2^{1/2}&v_2u_1^{-1/2}\\
0&0&u_1^{-1/2}\endmatrix \right)\cdot \left(\matrix \alpha_2^{-2}&0&0\\
0&\alpha_2&0\\
0&0&\alpha_2\endmatrix \right)\cdot\left(\matrix 1&t_2&x_{2}\\
0&1&0\\
0&0&1\endmatrix \right).$$
Let $y_i=\alpha_i^6$, $i=1,2$ then a Haar measure on $G/K$ may be given 
in terms of Langlands coordinates as follows
$$d\mu={{dy_1}\over {y_1^2}}\,{{dz_1}\over {u_1^2}}\,dx_1\,dt_1
={{dy_2}\over {y_2^2}}\,{{dz_2}\over {u_2^2}}\,dx_2\,dt_2$$ where
$z_1=v_1+iu_1$ and $z_2=v_2+iu_2$.

Let $\Gamma=\text{SL}(3,{\Bbb Z})$ acting on $G/K$, and ${\Cal  D}$ be a
 fundamental domain for $\Gamma$. Then by the theory of Eisenstein series,
$$L^2(\Gamma\backslash G/K)={\Bbb H}_0\oplus\Theta_0^{(1)}\oplus 
\Theta_0^{(2)}\oplus \Theta_{1,2}^{(2)}$$ where $H_0$ denotes the cusp 
forms of $\Gamma$, while
the Theta's may be defined as follows using Eisenstein series:

Associated to minimal parabolic subgroup $P_0$ we have the Eisenstein series 
$$E^0(Y;s,t):=\sum_{\gamma\in P_0\cap\Gamma\backslash\Gamma}
y_1(\gamma Y)^su_1(\gamma Y)^t.\eqno(7)$$ It  is known that this series
converges when $3\text{Re}(s)-\text{Re}(t)>2,\text{Re}(t)>1$ and admits a 
meromorphic to the whole $(s,t)$-space. Despite that 
there are many poles, but these which are of some interests to us are on 
the lines $t=1, 3s-t=2, 3s+t=3$. The residues at 
these poles are meromorphically continued Eisenstein series of one 
variable and generate the closed subspace $\Theta^{(1)}$.
One checks that $\Theta_0^{(2)}$ is simply the span of 
$E^0(Y,1/2+ir_1,1/2+ir_2).$ (A general fact due to Langlands.)

Now, let $\phi$ be an even cusp forms for $SL(2,{\Bbb Z})$ on the upper 
half-plane. Set 
$$E_i(Y;\phi;s):=\sum_{P_i\cap\Gamma\backslash \Gamma}
y_i(\gamma Y)^s\cdot \phi(z_i(\gamma Y)),\qquad i=1,2.\eqno(8)$$ 
These  series converge for $\text{Re}(s)>1$ and have meromorphic extensions 
on the whole
$s$-plane which have no poles on the line $(1/2,1]$. One checks that 
the space $\Theta_{1,2}$ generated by $E_i(Y;\phi;s),i=1,2$ for all $\phi$
coincides with $\Theta_{1,2}^{(2)}$, the closed space  spanned by $E_i$ 
along the line $\text{Re}(s)=1/2$. Indeed, one may also have a refined 
orthogonal decomposition of $\Theta_{1,2}^{(2)}$ according to that of $\phi$. 
For details, see [V] or Ch. III below.
\vskip 0.30cm
\noindent
{\bf II.3.2. Functional Equations}

There may be different ways to write down the functional equations. For us, to
indicate how delicate they are, we first introduce what we call the completed
Eisenstein series, then give  very elegent functional equations for them. 
(Please compare our approach to others, e.g, [V].)

Define the completed Eisenstein series associated to $E^0(Y;s,t)$ above
by $$\hat E^0(Y;s,t):=\xi(2t)\xi(3s-t)\xi(3s+t-1)\cdot E^0(Y;s,t).$$
Then one checks that we have the following very beautiful
\vskip 0.30cm
\noindent
{\bf Functional Equations.} {\it With the same notation as above,

\noindent
(i) $\hat E^0(Y;s,1-t)=\hat E^0(Y;s,t);$

\noindent
(ii) $\hat E^0(Y;{1\over 2}(1-s+t),{1\over 2}(-1+3s+t))=\hat E^0(Y;s,t);$

\noindent
(iii) $\hat E^0(Y;{1\over 2}(2-s-t),{1\over 2}(2-3s+t))=\hat E^0(Y;s,t);$

\noindent
(iv) $\hat E^0(Y;{1\over 2}(1-s+t),{1\over 2}(3-3s-t))=\hat E^0(Y;s,t);$

\noindent
(v) $\hat E^0(Y;{1\over 2}(2-s-t),{1\over 2}(3s-t))=\hat E^0(Y;s,t).$}
\eject
\vskip 0.30cm
\noindent
{\bf II.3.3. Fourier Expansions}
\vskip 0.30cm
To go further, we need to understand the Fourier expansion of Eisenstein 
series near cusps. However to motivate our discussion of non-abelian 
$L$-functions in Ch. III,  we concentrate on
 the Eisenstein series  $\hat E^0(Y;s,t)$ of the highest level, whose
residue gives the Epstein zeta function. 

Let us start with the simplest terms, i.e., the so-called constant terms 
appeared in the Fourier expansion for the cusps. By definition, for 
any measurable, locally $L^1$ function $f(Y)$ on $N_j$, set the constant 
term of $f$ along 
$P_j$ to be
$$f_{P_j}(Y):=\int_{\Gamma\cap N_j\backslash N_j}f(nY)dn,\qquad j=0,1,2.$$

\noindent
{\bf Proposition.} (See e.g, Venkov[V]) {\it With the same notation as above,} 
$$\align \hat E^0_{P_0}(Y;s,t)=&\xi(2t)\xi(3s-t)\xi(3s+t-1)\cdot y_1^su_1^t\\
&+\xi(2t)\xi(3s-t-1)\xi(3s+t-1)\cdot y_1^su_1^{1-t}\\
&+\xi(2t-1)\xi(3s-t-1)\xi(3s+t-2)\cdot
y_1^{{1\over 2}(1-s-t)}u_1^{{1\over 2}(2-3s+t)}\\
&+\xi(2t)\xi(3s-t-1)\xi(3s+t-1)\cdot y_1^{{1\over 2}(1-s-t)}u_1^{{1\over 2}
(3-3s-t)}\\
&+\xi(2t-1)\xi(3s-t)\xi(3s+t-2)\cdot
y_1^{{1\over 2}(2-s-t)}u_1^{{1\over 2}(3s-t)};\\
\hat E_{P_i}^0(Y;s,t)=&\xi(3s-t)\xi(3s+t-1)\cdot y_i^s\cdot \hat E(z_i,t)\\
&+\xi(2t)\xi(3s-t-1)\cdot y_i^{{1\over 2}(1-s-t)}\cdot \hat E(z_i,
{1\over 2}(3s+t-1))\\
&\xi(2t-1)\xi(3s-t)\cdot y_i^{{1\over 2}(2-s-t)}\cdot \hat E(z_i,{1\over 2}
(3-3s-t));\qquad i=1,2\endalign$$

For the proof, see, e.g., that of Lemmas 2 and 8 of [V]. (Note that here
we systematically use the completed Eisenstein series.)
\vskip 0.30cm
Next, let us recall the Fourier expansions of $\hat E^0(Y;s,t)$ along 
the parabolic subgroups $P_1$ and $P_2$. 
(In theory, we should also know the Fourier expansion along $P_0$. 
However, as the later calculation shows, with an induction on the rank, 
to see the non-abelian 
contributions, such detailed information is not needed: 
the terms involved will finally lead to a 
combination of abelian zeta functions just as what happens for
$SL_2$  in Chapter II.) We will follow [IT].
For this,  view rank 3 lattices of volume one as positive quadratic forms 
of determinant 1,
write $$Y=\left(\matrix U&0\\ 0&w\endmatrix \right)\left[\matrix I_2&x\\
0&1\endmatrix \right]=\left(\matrix I_2&0\\
x^t&1\endmatrix \right)\cdot\left(\matrix U&0\\ 0&w\endmatrix \right)
\cdot\left(\matrix I_2&x\\
0&1\endmatrix \right)$$ and define the first type of matrix k-Bessel 
function to be
$$k_{2,1}(Y;s_1,s_2;A):=\int_{X\in {\Bbb R}^{2\times 1}}
p_{-s_1,-s_2}\left(Y^{-1}\left[\matrix 1&0\\ x^t&I_2\endmatrix \right]\right)
\exp\Big(2\pi i {\text {Tr}}(A^t\cdot X)\Big)dX$$ 
for $(s_1,s_2)\in {\Bbb C}^2$, $Y\in {\Cal SP}_3, A\in {\Bbb R}^{2\times 1}$ 
and $p_{s_1,s_2}(Y):=|Y_1|^{s_1}|Y_2|^{s_2}$
where $Y_j\in {\Cal  SP}_j$ is the $j\times j$ upper left hand corner in 
$Y,\ j=1,2.$ Here as usual, we denote by ${\Cal  SP}_n$
the collection of rank $n$ positive quadratic forms of determinant 1.
  Set also 
$$\eqalign{\Lambda(s,r)=&\pi^{-(s-{r\over 2})}\Gamma(s-{r\over 2})
\pi^{-(s-{{1-r}\over 2})}\Gamma(s-{{1-r}\over 2}),\cr
\alpha_0=&{{\Lambda(s,r)}\over{B({1\over 2},{1\over 2}-r)}},\qquad \alpha_0'
={{\Lambda(s,r)}\over
{B({1\over 2},r-{1\over 2})}},\qquad
\alpha_{k\not=0}=\Lambda(s,r){{\sigma_{1-2r}(k)}\over{\zeta(2r)}},\cr
c(s,r)=&\xi(2r)\xi(2s-r)\xi(2s-1+r)\cdot E({U\over {\sqrt{|U|}}};r)|U|^{-s}
\qquad \text{with}\cr
 E(V;r)=&{1\over 2}\sum_{\text{gcd}(a)=1}V[a]^{-r},\qquad 
\text{Re}(r)>1.\cr}$$
\vskip 0.30cm
\noindent
{\bf Proposition.} ([IT]) {\it With the same notation as above, we have 
$$\align ~&\Lambda(s,r)\cdot E^0
\Big(\left(\matrix U&0\\ 0&w\endmatrix \right)
\left[\matrix I_2&x\\ 0&1\endmatrix \right];r,s\Big)\\
=&c(s,r)+c({{6-2s-3r}\over 4},s-{r\over 2})
+c({{3+3r-2s}\over 4},s-{{1-r}\over 2})\\
&+\sum_{A\in \text{SL}(2,{\Bbb Z})/P(1,1)}
\Big(\sum_{c, d_2\in {\Bbb Z}_{>0},d_1\in {\Bbb Z}\backslash \{0\}}
\Big[\alpha_0'\,c^{2-2s-r}\,d_2^{r-2s}
\exp\Big(2\pi ix^t A\cdot{cd_1\choose 0} \Big)\\
&\hskip 4.0cm 
\cdot k_{2,1}(\left(\matrix A^{-1}UA^{-t}&0\\ 0&w\endmatrix \right);
s-{r\over 2},r;\pi {cd_1\choose 0})\\
&\quad+\alpha_0\,c^{1-2s+r}\,d_2^{1-r-2s}
\exp\Big(2\pi ix^t A{cd_1\choose 0} \Big)\\
&\hskip 4.0cm  
\cdot k_{2,1}
\Big(\left(\matrix A^{-1}UA^{-t}&0\\ 0&w\endmatrix \right);
s-{{1-r}\over 2},1-r;\pi{cd_1\choose 0}\Big)\Big]\\
&+\sum_{k\not=0}\sum_{c, d_2\in {\Bbb Z}_{>0},d_2|k,d_1\in 
{\Bbb Z}\backslash \{0\}}\alpha_k\,c^{2-2s-r}\,d_2^{r-2s}
\exp\Big(2\pi ix^t A{cd_1\choose {{ck}/{d_2}}} \Big)\\
&\hskip 4.0cm 
\cdot k_{2,1}\Big(\left(\matrix A^{-1}UA^{-t}&0\\ 0&w\endmatrix \right);
s-{r\over 2},r;\pi {cd_1\choose {{ck}/{d_2}}}\Big)\Big),\endalign$$ 
where $P(1,1)$ is the subgroup of upper triangle matrices of 
determinant 1.
Similar Fourier expansion holds for $E^0(E(z;s),t)$ with respect to $P_2$.}
\vskip 0.30cm
\noindent
{\bf II.3.4. Non-abelian Contributions}
\vskip 0.30cm
To see where non-abelian contribution come, we go as follows. 
(We reminder the reader that the purpose of this example is to motivate
the discussion of non-abelian $L$-functions in Ch. III.)
 
First, for simplicity, consider the geometric truncated fundamental domain of 
$\Gamma:={SL}(3,{\Bbb Z})$ obtained by cutting off
the cusp regions corresponding to $P_1,P_2$ and $P_0$.
Put $\Gamma_j=\Gamma\cap P_j, j=0,1,2$ and $\Gamma_{N_0}=\Gamma\cap N_0$. 
Then the fundamental domain $F_*$ in
$S:={SL}(3,{\Bbb R})/{SO}(3,{\Bbb R})$ for the groups 
$*=\Gamma_0,\Gamma_1,\Gamma_2,\Gamma_{N_0}$ may be choosen to be
$$\eqalign{F_{N_0}:=&\{Y\in S:y_1>0,u_1>0,-1/2<v_1,x_1,t_1<1/2\};\cr
 F_0:=&\{Y\in F_{N_0}:v_1+x_1>0,v_1+t_1>0,x_1+t_1>0\};\cr
F_j:=&\{Y\in F_0:v_j^2+u_j^2\geq 1\},\qquad j=1,2.\cr} $$ 
With this, one checks that there exists a compact set $F^0\subset S$ such that
$$F_1\cap F_2=F^0\cup F$$ where $F$ denotes the fundamental domain of 
$\Gamma$. That is to say, the cusp regions
for $P_j$, $j=0,1,2$ in the fundamental region $F$ of ${SL}(3,{\Bbb Z})$ 
may be read from $F_1$ and $F_2$.

Set also $D_j^T:=\{Y=\{y_j;z_j,x_j,t_j\}\in F:y_j\geq T\}, j=1,2$ for 
sufficient large $T>0$. With this, then we get a geometric truncated 
compact subset in $F$ by cutting off the neighborhood of cusps along 
$F_1$ and $F_2$, so as to get
$F_T:=F\backslash (D_1^T \cup D_2^T)$. Note that $D_0^T:=D_1^T\cap D_2^T$ 
gives a neighborhood for the cusps with respect to
$P_0$. Thus, we may analytically understand this geometric truncation as
$$1_{F_T}=1_F-1_{D_1^T}-1_{D_2^T}+1_{D_0^T},$$ which is compactible 
with the truncations in the next chapter.
\vskip 0.30cm
Secondly, let us simply look at the contributions of standard parabolic 
subgroups so as to get the analytic truncation
$$\eqalign{&\Lambda_T\hat E^0(Y;s,t):\cr
=&\hat E^0(Y;s,t)-\hat E_{P_1}^0(Y;s,t)\cdot 1_{D_1^T}
-\hat E_{P_2}^0(Y;s,t)\cdot 1_{D_2^T}+\hat E_{P_0}^0(Y;s,t)\cdot 1_{D_0^T}\cr
=&\Big(\big(\hat E^0(Y;s,t)-\hat E_{P_1}^0(Y;s,t)\cdot 1_{D_1^T}\big)
+\big(\hat E^0(Y;s,t)-\hat E_{P_2}^0(Y;s,t)\cdot 1_{D_2^T}\big)\Big)\cr
&\qquad-\Big(\hat E^0(Y;s,t)-\hat E_{P_0}^0(Y;s,t)\cdot 1_{D_0^T}\Big)\cr
=&H_{P_1}^0(Y;s,t)+H_{P_2}^0(Y;s,t)-H_{P_0}^0(Y;s,t).\cr}$$ Here
$$H_{P_j}^0(Y;s,t):=\hat E^0(Y;s,t)-\hat E_{P_j}^0(Y;s,t)\cdot 1_{D_j^T},\qquad j=1,2,0$$ denotes the non-constant part of the corresponding Fourier 
expansion.

Thirdly, we want to know the integration $\int_{F_T}\hat E^0(Y;s,t)d\mu(Y)$. For this, we go as follows:
$$\eqalign{\int_{F_T}\hat E^0(Y;s,t)d\mu(Y)
=&\int_{F}\Lambda^T\hat E^0(Y;s,t)d\mu(Y)
-\int_{F\backslash F_T}\Lambda^TE^0(Y;s,t)d\mu(Y)\cr
=&\int_{F}\Lambda^T\hat E^0(Y;s,t)d\mu(Y)
-\int_{F^T}\Lambda^T\hat E^0(Y;s,t)d\mu(Y),\cr} $$ where 
$F^T:=F\backslash F_T=D^T_1\cup D^T_2$.
\vskip 0.30cm
Finally, let us look at the structure of this latest expression:

\noindent
(A) ({\it Abelian Part:  Application of Rankin-Selberg Method}) 
By the Rankin-Selberg method, in particular, the version generalized 
by Langlands and Arthur, (see Ch. IV below where a detailed formula is given,)
the part $\int_{F}\Lambda^T\hat E^0(Y;s,t)d\mu(Y)$,   
being the integration of analytic
truncated Eisenstein series on the whole fundamental domain of 
${SL}(3,{\Bbb Z})$,  is essentially abelian;
\vskip 0.30cm
Thus, it suffices to know the structure of  
$\int_{F^T}\Lambda^T\hat E^0(Y;s,t)d\mu(Y)$. Clearly,
$$\eqalign{&\int_{F^T}\Lambda^T\hat E^0(Y;s,t)d\mu(Y)\cr
=&\int_{D_1^T}\Lambda^T\hat E^0(Y;s,t)d\mu(Y)
+\int_{D_2^T}\Lambda^T\hat E^0(Y;s,t)d\mu(Y)
-\int_{D_0^T}\Lambda^T\hat E^0(Y;s,t)d\mu(Y)\cr
=&\int_{D_1^T}\Big(H_{P_1}^0(Y;s,t)+H_{P_2}^0(Y;s,t)
-H_{P_0}^0(Y;s,t)\Big)d\mu(Y)\cr
&\qquad+\int_{D_2^T}\Big(H_{P_1}^0(Y;s,t)+H_{P_2}^0(Y;s,t)
-H_{P_0}^0(Y;s,t)\Big)d\mu(Y)\cr
&\qquad\qquad-\int_{D_0^T}\Big(H_{P_1}^0(Y;s,t)+H_{P_2}^0(Y;s,t)
-H_{P_0}^0(Y;s,t)\Big)d\mu(Y)\cr
=&I_1^T(s,t)+I_2^T(s,t)-I_0^T(s,t),\cr} $$ where
$$I_j^T(s,t):=\int_{D_j^T}\Big(H_{P_1}^0(Y;s,t)+H_{P_2}^0(Y;s,t)
-H_{P_0}^0(Y;s,t)\Big)d\mu(Y),j=0,1,2.$$
\noindent
(B) ({\it Terms obtained from Lower Rank Non-Abelian Zeta:  
Induction on the Rank}) 

Consider the integrations 
$$I_j^T(s,t):=\int_{D_1^T}\Big(H_{P_1}^0(Y;s,t)+H_{P_2}^0(Y;s,t)
-H_{P_0}^0(Y;s,t)\Big)d\mu(Y),j=0,1,2.$$
If  the fundamental domain $F$ is chosen so that $F$ is 
of exact box shape as $Y$ approaches to 
all levels of cusps, we have $$\int_{D_j^T}H_{P_j}(Y;s,t)=0.$$ 
(This choice of fundamental domain is possible by a result of 
Grenier recalled in [T]. From now on, we always assume 
that this condition for
the fundamental domain holds.) Then what left is to consider
the following integrations:
$$\eqalign{{II}_1^T(s,t):=&\int_{D_1^T}\Big(H_{P_2}^0(Y;s,t)
-H_{P_0}^0(Y;s,t)\Big)d\mu(Y);\cr
{II}_2^T(s,t):=&\int_{D_2^T}\Big(H_{P_1}^0(Y;s,t)-
H_{P_0}^0(Y;s,t)\Big)d\mu(Y);\cr
{III}^T(s,t):=&\int_{D_0^T}\Big(H_{P_1}^0(Y;s,t)+H_{P_2}^0(Y;s,t)\Big)d\mu(Y).
\cr}$$

By the structures of $D_1^T, D_2^T$, we then see that ${II}_i^T(s,t)$ are 
essentially  rank two zeta functions, which may be understood via an 
induction argument.
So we are left with only 
$${III}^T(s,t):=\int_{D_0^T}\Big(H_{P_1}^0(Y;s,t)+
H_{P_2}^0(Y;s,t)\Big)d\mu(Y),$$ which in the case of rank 3, 
is the only essential non-abelian contribution, as indicated below.
\vskip 0.30cm
\noindent
(C) ({\it Essential Non-abelain Contributions:  New Ingredients}) 

${III}^T(s,t)$ gives naturally the non-abelian contribution:
Indeed, by definition, $H_{P_i}^0(Y;s,y), i=1,2$ are the  non-constant 
terms of the 
Fourier expansion of $\hat E^0(Y;s,t)$ along $P_i$, $i=1,2$.
Moreover, the integration is taken to be the cusp region corresponding to the
minimal parabolic subgroup $P_0$, which are proper subgroups of $P_1$ and 
$P_2$.  By the result of Imai and Terras 
cited above, they consist of matrix version of
k-Bessal functions. 

Up to this point, we see then clearly what a new direction should we go:
while the classical Rankin-Selberg method gives a natural way to single out
abelian contributions out of non-abelian data, we should treat both 
abelian and non-abelian parts uniformly in order to see their intrinsic 
structures.
\vskip 0.5cm
\centerline {\li Chapter III. New Non-Abelian L-Functions}
\vskip 0.30cm
\centerline
{\li III.1. Geo-Arithmetical Truncation}
\vskip 0.30cm
\noindent
{\bf III.1.1. Slopes, Canonical Filtrations}

Following Lafforgue [Laf], we call an abelian category ${\Cal A}$ together with 
two additive
morphisms $$\text {rank}:{\Cal A}\to {\Bbb N},
\qquad \text{deg}:{\Cal A}\to {\Bbb R}$$ a category with slope structure. 
In particular, for non-zero $A\in {\Cal A}$, 

\noindent
(1) define the slope of $A$ by 
$$\mu(A):={{\text{deg}(A)}\over {\text{rank}A}};$$
  
\noindent
(2) If $0=A_0\subset A_1\subset\cdots\subset A_l=A$ is a filtration of $A$ in 
${\Cal A}$ with $\text {rank}(A_0)<\text {rank}(A_1)<\cdots<
\text {rank}(A_l)$, then define the associated polygon to be the continuous 
function $[0,\text {rank}A]\to {\Bbb R}$ such that

\noindent
(i) its values at 0 and $\text {rank}(A)$ are 0;

\noindent
(ii)  it is affine with slope $\mu(A_i/A_{i-1})-\mu(A)$ on the intervals 
$[\text {rank}(A_{i-1}),\text {rank}(A_i)]$ for all $1\leq i\leq l$,

\noindent
(3) If $\frak a$ is a collection of subobjects of $A$ in ${\Cal A}$, then 
$\frak a$ is said to be nice if

\noindent
(i) $\frak a$ is stable under intersection and finite summation;

\noindent
(ii) $\frak a$ is Noetherian in the sense that every increasing chain 
of elements
in $\frak a$ has a maximal element in $\frak a$;

\noindent
(iii) if $A_1\in\frak a$ then $A_1\not=0$ if and only if 
$\text {rank}(A_1)\not=0$;

\noindent
(iv) for $A_1,A_2\in\frak a$ with $\text {rank}(A_1)=\text {rank}(A_2)$. Then
$A_1\subset A_2$ is proper implies that $\text {deg}(A_1)<\text {deg}(A_2);$

\noindent
(4) For any nice $\frak a$, set 
$\mu^+(A):=\text {sup}\{\mu(A_1):A_1\in\frak a,
\text {rank}(A_1)\geq 1\}$,
$\mu^-(A):=\text {inf}\{\mu(A/A_1):A_1\in\frak a,
\text {rank}(A_1)< \text {rank}(A) \}$. Then we say $(A,\frak a)$ is 
semi-stable if $\mu^+(A)=\mu(A)=\mu^-(A)$. Moreover if $\text {rank}(A)=0$ 
set also $\mu^+(A)=-\infty$ and $\mu^-(A)=+\infty$.
\vskip 0.30cm
\noindent
{\bf Proposition 1.} ([Laf]) {\it Let ${\Cal A}$ be a 
category with slope, $A$ an object in ${\Cal A}$ and $\frak a$ a nice 
family of subobjects of $A$ in ${\Cal A}$. Then

\noindent
(1) $A$ admits a unique filtration $0=
\bar A_0\subset \bar A_1\subset\cdots\subset\bar A_l=A$
with elements in $\frak a$  such that

\noindent
(i) $\bar A_i,0\leq i\leq k$ are maximal in $\frak a$;

\noindent
(ii) $\bar A_i/\bar A_{i-1}$ are semi-stable;

\noindent
(iii) $\mu(\bar A_1/\bar A_{0})>\mu(\bar A_2/\bar A_{1}>\cdots>
\mu(\bar A_k/\bar A_{k-1})$.

\noindent
We call such a filtration the canonical filtration of $(A,\frak a)$;

\noindent
(2) All polygons of filtrations of $A$ with elements 
in $\frak a$ are bounded from above by $\bar p$, where 
 $\bar p:=\bar p^A$ be the associated polygon for the canonical 
filtration in (1);

\noindent
(3) For any $A_1\in \frak a, \text {rank}(A_1)\geq 1$ implies
$\mu(A_1)\leq \mu(A)+\bar p(\text {rank}(A_1))/\text {rank}(A_1);$

\noindent
(4) The polygon $\bar p$ is convex with maximal slope $\mu^+(A)-\mu(A)$ 
and minimal slope $\mu^-(A)-\mu(A)$;

\noindent
(5) If $(A',\frak a')$ is another pair, and $u:A\to A'$ is a homomorphism 
such that $\text {Ker}u\in\frak a$ and $\text{Im}u\in\frak a'$. Then
$\mu^-(A)\geq\mu^+(A')$ implies that $u=0$.}
\vskip 0.30cm
\noindent
{\it Proof.} This results from a Harder-Narasimhan type filtration 
consideration. A detailed proof may be found at pp. 87-88 in [Laf].  
\vskip 0.30cm
\noindent
{\bf Proposition 2.} {\it  Let $F$ be a number field. Then

\noindent
(1) the abelian category of hermitian vector sheaves on 
$\text {Spec}{\Cal O}_F$ together with the natural rank and the 
Arakelov degree is a category with slopes;

\noindent
(2) For any hermitian vector sheaf $(E,\rho)$, $\frak a$ consisting 
of pairs $(E_1,\rho_1)$ with $E_1$ sub vector sheaves of $E$ and 
$\rho_1$ the restriction of $\rho$ forms a nice family.}
\vskip 0.30cm
\noindent
{\it Proof.} (1) is obvious, while (2) is a direct consequence of 
the following standard fact: For a fixed $(E,\rho)$,
$\{\text {deg}(E_1,\rho_1):(E_1,\rho_1)\in \frak a\}$ is discrete 
subset of ${\Bbb R}$.
\vskip 0.30cm
Thus in particular we  get the canonical filtration  of Harder-Narasimhan
 type for 
hermitian vector sheaves over $\text {Spec}{\Cal O}_F$, which are indeed 
${\Cal O}_F$-lattices in $({\Bbb R}^{r_1}\times 
{\Bbb C}^{r_2})^{r=\text{rank}(E)}$ induced from the natural embedding
$F^r\hookrightarrow ({\Bbb R}^{r_1}\times {\Bbb C}^{r_2})^r$ where $r_1$ 
(resp. $r_2$) denotes the real (resp. complex) embeddings of $F$.
\vskip 0.30cm
\noindent
{\bf III.1.2. Canonical Polygons and Geo-Arithmetical Truncation}

Let $X=\text {Spec}{\Cal O}_F$. If $E$ is a vector sheaf of rank $r$ over $X$,
 i.e, a locally free ${\Cal O}_F$-sheaf of rank $r$, denote by $E_F$ the fiber 
of $E$ at the generic point $\text {Spec}(F)\hookrightarrow 
\text {Spec}{\Cal O}_F$ of $X$ ($E_F$ is an $F$-vector space of dimension $r$), 
and for each $v\in S_f$, set $E_{{\Cal O}_v}:=
H^0(\text {Spec}{\Cal O}_{F_v},E)$ a free ${\Cal O}_v$-module of rank 
$r$. In particular, we have a canonical isomorphism:
$$\text{can}_v:F_v\otimes_{{\Cal O}_v}E_{{\Cal O}_v}\simeq F_v\otimes_FE_F.$$
 
Now if $E$ is a vector sheaf of rank $r$ over $X$ equipped with a basis 
$\alpha_F:F^r\simeq E_F$ of its generic fiber and a basis 
$\alpha_{{\Cal O}_v}:{\Cal O}_v^r\simeq E_{{\Cal O}_v}$ for any 
$v\in S_f:=S_{\text{fin}}$, the elements $g_v:=(F_v\otimes_F\alpha_F)^{-1}\circ\text{can}_v
\circ (F_v\otimes_{{\Cal O}_v}\alpha_{{\Cal O}_v})\in GL_r(F_v)$ for all 
$v\in S_f$ define an element $g_{\Bbb A}:=(g_v)_{v\in S_f}$ of 
$GL_r({\Bbb A}_f)$, i.e, for almost every $v$ we have $g_v\in 
GL_r({\Cal O}_v)$. By this construction, we obtain a bijection from 
the set of isomorphism classes of triples $(E;\alpha_F;
(\alpha_{{\Cal O}_v})_{v\in S_f})$ as above onto $GL_r({\Bbb A}_f)$. 
Moreover, if $r\in GL_r(F), k\in GL_r({\Cal O}_F)$ and if this bijection
 maps the triple $(E;\alpha_F;(\alpha_{{\Cal O}_v})_{v\in S_f})$ onto 
$g_{\Bbb A}$, the same map maps the triple $(E;\alpha_F\circ r^{-1};
(\alpha_{{\Cal O}_v}\circ k_v)_{v\in S_f})$ onto $rg_{\Bbb A}k$. 
Therefore the above bijection induces a bijection between the set of 
isomorphism classes of vector sheaves of rank $r$ on 
$\text {Spec}{\Cal O}_F$ and the double coset space
$GL_r(F)\backslash GL_r({\Bbb A}_f)/GL_r({\Cal O}_F)$.

More generally, let $r=r_1+\cdots+r_s$ be a partition $I=(r_1,\cdots,r_s)$ of 
$r$ and let $P_I$ be the corresponding standard parabolic subgroup of $GL_r$. 
Then we have a natural bijection from the set of isomorphism classes of 
triple $(E_*;\alpha_{*,F}:(\alpha_{*,{\Cal O_v}})_{v\in S_f})$ onto 
$P_I({\Bbb A}_f)$, where $E_*:=\big((0)=E_0\subset E_1\subset\cdots\subset
 E_s\big)$ is a filtration of vector sheaves of rank $(r_1,r_1+r_2,\cdots, 
r_1+r_2+\cdots+r_s=r)$ over $X$, (i.e, each $E_j$ is a vector sheaf of rank 
$r_1+r_2+\cdots+r_j$ over $X$ and each quotient $E_j/E_{j-1}$ is torsion 
free,) which is equipped with an isomorphism of filtrations of $F$-vector 
spaces
$$\alpha_{*,F}:\big((0)=F_0\subset F^{r_1}\subset\cdots\subset F^{r_1+r_2+
\cdots+r_s=r}\big)\simeq (E_*)_F,$$ and with an isomorphism of filtrations 
of free ${\Cal O}_v$-modules
$$\alpha_{*,{\Cal O}_v}:\big((0)\subset {\Cal O}_v^{r_1}\subset\cdots\subset 
{\Cal O}_v^{r_1+r_2+\cdots+r_s=r}\big)\simeq (E_*)_{{\Cal O}_v},$$ for every 
$v\in S_f$. Moreover this bijection induces a bijection between the set of 
isomorphism classes of the filtrations of vector sheaves of rank 
$(r_1,r_1+r_2,\cdots, r_1+r_2+\cdots+r_s=r)$ over $X$ and the double coset 
space $P_I(F)\backslash P_I({\Bbb A}_f)/P_I({\Cal O}_F)$. The natural 
embedding $P_I({\Bbb A}_f)\hookrightarrow P_I({\Bbb A})$ (resp. the 
canonical projection $P_I({\Bbb A}_f)\to M_I({\Bbb A}_f)\to GL_{r_j}
({\Bbb A}_f)$ for $j=1,\cdots,s$, where $M_I$ denotes the standard Levi 
of $P_I$) admits the modular interpretation
$$(E_*;\alpha_{*,F}:(\alpha_{*,{\Cal O_v}})_{v\in S_f})\mapsto 
(E_s;\alpha_{s,F}:(\alpha_{s,{\Cal O_v}})_{v\in S_f})$$ (resp. 
$$(E_*;\alpha_{*,F}:(\alpha_{*,{\Cal O_v}})_{v\in S_f})\mapsto 
(\text{gr}_j(E_*);\text{gr}_j(\alpha_{*,F}), \text{gr}_j
(\alpha_{*,{\Cal O_v}})_{v\in S_f}),$$ where $\text{gr}_j(E_*):=
E_j/E_{j-1}$,
$\text{gr}_j(\alpha_{*,F}):F^{r_j}\simeq \text{gr}_j(E_*)_F$ and 
$\text{gr}_j(\alpha_{*,{\Cal O_v}}):{\Cal O_v}^{r_j}\simeq \text{gr}_j
(E_*)_{\Cal O_v}$, $v\in S_f$ are induced by $\alpha_{*,F}$ and 
$\alpha_{*,{\Cal O_v}}$ respectively.)

Moreover, any $g=(g_f;g_\infty)\in GL_r({\Bbb A}_f)\times 
GL_r({\Bbb A}_\infty)=GL_r({\Bbb A})$ gives first a rank $r$ 
vector sheaf $E_g$ on $\text{Spec}{\Cal O}_F$, which via the 
embedding $E_F\hookrightarrow ({\Bbb R}^{r_1}\times {\Bbb C}^{r_2})^r$ 
gives a discrete subgroup, a  free rank $r$ ${\Cal O}_F$-module. In particular, 
$g_\infty=(g_\sigma)$ then induces a natural metric  on $E_g$ by twisting 
the standard one on $({\Bbb R}^{r_1}\times {\Bbb C}^{r_2})^r$ via the 
linear transformation induced from $g_\infty$. As a direct consequence, 
see e.g, [L3], $$\text{deg}(E_g,\rho_g)=-\log\big(N(\text{det}g)\big)$$ 
with $N:GL_1({\Bbb A}_F)={\Bbb I}_F\to {\Bbb R}_{>0}$ the standard norm 
map of the idelic group of $F$.

With this, for $g=(g_f;g_\infty)\in GL_r({\Bbb A})$ and a parabolic subgroup 
$Q$ of $GL_r$, denote by $E_*^{g;Q}$ the filtration of the vector sheaf 
$E_{g_f}$ induced by the parabolic subgroup $Q$. Then we have a filtration 
of hermitian vector sheaves $(E_*^{g;Q},\rho_*^{g;Q})$ with the hermitian 
metrics $\rho_j^{g;Q}$
on $E_j^{g;Q}$ obtained via the restrictions of $\rho_{g_\infty}$.

Now introduce an associated polygon $p_Q^g:[0,r]\to {\Bbb R}$ by the following 3 conditions:

\noindent
(i) $p_Q^g(0)=p_Q^g(r)=0$;

\noindent
(ii) $p_Q^g$ is affine on the interval 
$[\text{rank}E_{i-1}^{g;Q},\text{rank}E_{i}^{g;Q}]$; and 

\noindent
(iii) for all indices $i$, 
$$p_Q^g(\text{rank}E_{i}^{g;Q})=\text{deg}(E_{i}^{g;Q},\rho_{i}^{g;Q})-
{{\text{rank}E_{i}^{g;Q}}\over r}\cdot \text{deg}(E_{g},\rho_{g}).$$
Then by Prop. 2 of III.1.1, there is a unique convex polygon $\bar p^g$ 
which bounds all $p_Q^g$ from above for all parabolic subgroups $Q$ for $GL_r$. 
Moreover there exists a parabolic subgroup $\bar Q^g$ such that 
$p_{{\bar Q}^g}^g=\bar p^g$.

Hence we have the following foundamental
\vskip 0.30cm
\noindent
{\bf Main Lemma.} {\it For any fixed polygon $p:[0,r]\to {\Bbb R}$ and any $d\in {\Bbb R}$, the subset
$$\{g\in GL_r(F)\backslash GL_r({\Bbb A}):\text{deg}\,g=d,\bar p^g\leq p\}$$ 
is compact.}
\vskip 0.30cm
\noindent
{\it Proof.} This is a restatement of the classical reduction theory. 
It consists of two parts: In terms of ${\Cal O}_F$-lattices, 
by fixing the degree, we get a fixed volume for the free rank $r$
${\Cal O}_F$-lattice corresponding to $(E_g,\rho_g)$, by the Arakelov 
Riemann-Roch [L3]. Thus the condition $\bar p^g\leq p$ gives an upper bound for
the volumes of all the sublattices of $(E_g,\rho_g)$ and hence all the 
Minkowski successive minimas.

On the other hand, the reduction theory in terms of adelic language also 
tells us that the fiber of the natural map from $GL_r(F)\backslash 
GL_r({\Bbb A})$ to isomorphism classes of ${\Cal O}_F$-lattices are all 
compact. (This is in fact generally true for all reductive groups, a 
result due to Borel [Bo1,2]. However, our case
where only $GL_r$ is involved is rather obvious: essentially, the fibers are
the compact subgroup 
$GL_r({\Cal O}_F)\times SO_r({\Bbb R})^{r_1}\times SU_r({\Bbb C})^{r_2}$.)
This completes the proof.
\vskip 0.30cm 
Similarly yet more generally, for a fixed parabolic subgroup $P$ of $GL_r$ 
and $g\in GL_r({\Bbb A})$, there is a unique maximal element $\bar p_P^g$ 
among all $p_Q^g$, where $Q$ runs over all parabolic subgroups of $GL_r$ 
which are contaiend in $P$.  And we have
\vskip 0.30cm
\noindent
{\bf Main Lemma}$'$. {\it For any fixed polygon $p:[0,r]\to {\Bbb R}$, $d\in {\Bbb R}$ 
and  any standard parabolic subgroup $P$ of $GL_r$, the subset
$$\{g\in GL_r(F)\backslash GL_r({\Bbb A}):\text{deg}\,g=d,\bar p^g_P\leq p, 
p^g_P\geq -p\}$$ is compact.}
\vskip 0.30cm
\noindent
{\bf III.1.3. Relation to Arthur's Analytic Truncation}

Let $p,q:[0,r]\to {\Bbb R}$ be two polygons and $P$ a standard parabolic 
subgroup of $GL_r$. Then we say $q>_Pp$ if for any $1\leq i\leq |P|$,
$$q(\text{rank}E_i^P)>p(\text{rank}E_i^P)$$ where $(r_1,\cdots,r_{|P|})$ 
denotes the partition of $r$ corresponding to $P$. Also as usual denote by 
${\bold 1}$ the characteristic function of the variable $g\in GL_r({\Bbb A})$. 
For example,
$${\bold 1}(\bar p^g\leq p)(g)=\cases 1,&\text {if}\ p^g\leq p\\
0&\text{otherwise}.\endcases$$
\vskip 0.30cm
\noindent
{\bf Proposition.} {\it For any convex polygon $p:[0,r]\to {\Bbb R}$, as a 
function of $g\in GL_r({\Bbb A})$,
$${\bold 1}(\bar p^g\leq p)=\sum_{P\supset P_0}(-1)^{|P|-1}
\sum_{\delta\in P(F)\backslash GL_r(F)}{\bold 1}(p_P^{\delta g}>_Pp).$$ 
Here $P$ runs over all standard parabolic subgpoups of $GL_r$.}

\noindent
{\it Proof.} For curves over finite fields, this result is due to Lafforgue, 
which parallel to Arthur's analytic truncation in trace formula. 

The basic idea is very simple: When the intersection stability breaks down, 
there 
should be naturally a parabolic subgroup which takes  responsibility.
With this, the detailed proof can easily be given 
 with the reference to pp.221-222 in [La]. We leave it to 
the reader.
\vskip 0.30cm
\centerline
{\li III.2. Non-Abelian L-Functions}
\vskip 0.30cm
\noindent
{\bf III.2.1. Choice of Moduli Spaces}

For the number field $F$ with discriminant $\Delta_F$, and for a fixed 
$r\in {\Bbb Z}_{>0}$, we take the moduli space to be 
$${\Cal M}_{F,r}^{\leq p}[|\Delta_F|^{r\over 2}]:=\{g\in GL_r(F)\backslash 
GL_r({\Bbb A}):\text{deg} g=-{r\over 2}\log|\Delta_F|,\bar p^g\leq p\}$$ 
for a fixed convex polygon $p:[0,r]\to{\Bbb R}$. Also 
we denote by $d\mu$ the induced Tamagawa measures on 
${\Cal M}_{F,r}^{\leq p}[|\Delta_F|^{r\over 2}]$. For example,  $p=0$ 
coincides the semi-stability introduced in Ch. I.

More generally, for any standard parabolic subgroup $P$ of $GL_r$, we 
introduce the moduli spaces $${\Cal M}_{F,r}^{P;\leq p}[|\Delta_F|^{r\over 2}]
:=\{g\in P(F)\backslash GL_r({\Bbb A}):\text{deg} g=-{r\over 2}\log|\Delta_F|,
\bar p_P^g\leq p, \bar p_P^g\geq -p\}.$$ By the discussion in III.1, these 
moduli spaces 
${\Cal M}_{F,r}^{P;\leq p}[|\Delta_F|^{r\over 2}]$ are all compact, a key
 property which plays a central role in our definition of non-abelian 
$L$-functions below.
\vskip 0.30cm
\noindent
{\bf III.2.2. Choice of Eisenstein Series: First Approach to Non-Abelian 
$L$-Function}

To faciliate our ensuing discussion, we start with  some preperations. We 
follow [MW] closely.

Fix a connected reduction group $G$ defined over $F$, denote by $Z_G$ its 
center. Fix a minimal parabolic subgroup $P_0$ of $G$. Then $P_0=M_0U_0$, 
where as usual we fix once and for all the Levi $M_0$ and  the unipotent
 radical $U_0$. A parabolic subgroup $P$ is $G$ is called standard if 
$P\supset P_0$. For such groups write $P=MU$ with $M_0\subset M$ the 
standard Levi and $U$ the unipotent radical. Denote by $\text {Rat}(M)$ 
the group of rational characters of $M$, i.e, the morphism $M\to {\Bbb G}_m$ 
where ${\Bbb G}_m$ denotes the multiplicative group. Set 
$$\frak a_M^*:=\text {Rat}(M)\otimes_{\Bbb Z}{\Bbb C},\qquad \frak a_M
:=\text{Hom}_{\Bbb Z}(\text {Rat}(M),{\Bbb C}),$$ and $$\text {Re}\frak 
a_M^*:=\text {Rat}(M)\otimes_{\Bbb Z}{\Bbb R},\qquad \text{Re}\frak a_M
:=\text{Hom}_{\Bbb Z}(\text {Rat}(M),{\Bbb R}).$$ For any $\chi\in 
\text {Rat}(M)$, we obtain a (real) character $|\chi|:M({\Bbb A})\to 
{\Bbb R}^*$ defined by $m=(m_v)\mapsto m^{|\chi|}:=\prod_{v\in S}
|m_v|_v^{\chi_v}$ with $|\cdot|_v$ the $v$-absolute values. Set then 
$M({\Bbb A})^1:=\cap_{\chi\in \text {Rat}(M)}\text{Ker}|\chi|$, which is 
a normal subgroup of $M({\Bbb A})$. Set $X_M$ to be the group of complex 
characters which are trivial on $M({\Bbb A})^1$. Denote by 
$H_M:=\log_M:M({\Bbb A})\to \frak a_M$ the map such that 
$\forall\chi\in \text {Rat}(M)\subset \frak a_M^*,\langle\chi,
\log_M(m)\rangle:=\log(m^{|\chi|})$. Clearly, $$M({\Bbb A})^1
=\text{Ker}(\log_M);\qquad \log_M(M({\Bbb A})/M({\Bbb A})^1)\simeq 
\text{Re}\frak a_M.$$ Hence in particular there is a natural isomorphism 
$\kappa:\frak a_M^*\simeq X_M.$
Set $$\text{Re}X_M:=\kappa (\text{Re}\frak a_M^*),\qquad \text{Im}X_M:=
\kappa (i\cdot \text{Re}\frak a_M^*).$$ Moreover define our working space 
$X_M^G$ to be the subgroup of $X_M$ consisting of complex characters of 
$M({\Bbb A})/M({\Bbb A})^1$ which are trivial on $Z_{G({\Bbb A})}$.

Fix a maximal compact sungroup ${\Bbb K}$ such that for all standard parabolic
subgroups $P=MU$ as above, $P({\Bbb A})\cap{\Bbb K}=M({\Bbb A})
\cap{\Bbb K}\cdot
U({\Bbb A})\cap{\Bbb K}.$ Hence we get the Langlands decomposition 
$G({\Bbb A})=M({\Bbb A})\cdot U({\Bbb A})\cdot {\Bbb K}$. Denote by 
$m_P:G({\Bbb A})\to M({\Bbb A})/M({\Bbb A})^1$ the map $g=m\cdot n\cdot
 k\mapsto M({\Bbb A})^1\cdot m$ where $g\in G({\Bbb A}), m\in M({\Bbb A}), 
n\in U({\Bbb A})$ and 
$k\in {\Bbb K}$.
 
Fix Haar measures on $M_0({\Bbb A}), U_0({\Bbb A}), {\Bbb K}$ respectively 
such that 

\noindent
(1) the induced measure on $M(F)$ is the counting measure and the volume of
 the induced measure on $M(F)\backslash M({\Bbb A})^1$ is 1. (Recall that 
it is  a fundamental fact that $M(F)\backslash M({\Bbb A})^1$ is compact.)

\noindent
(2) the induced measure on $U_0(F)$ is the counting measure and the volume 
of $U(F)\backslash U_0({\Bbb A})$ is 1. (Recall that being unipotent radical, 
$U(F)\backslash U_0({\Bbb A})$ is compact.)

\noindent
(3) the volume of ${\Bbb K}$ is 1.

Such measures then also induce Haar measures via $\log_M$ to $\frak a_{M_0}, 
\frak a_{M_0}^*$, etc. Furthermore, if we denote by $\rho_0$ the half of the 
sum of the positive roots of  the maximal split torus $T_0$ of the central 
$Z_{M_0}$
of $M_0$, then $$f\mapsto \int_{M_0({\Bbb A})\cdot U_0({\Bbb A})\cdot 
{\Bbb K}}f(mnk)\,dk\,dn\,m^{-2\rho_0}dm$$ defined for continuous functions 
with compact supports on $G({\Bbb A})$  defines a Haar measure $dg$ on 
$G({\Bbb A})$. This in turn gives measures on $M({\Bbb A}), U({\Bbb A})$ 
and hence on $\frak a_{M}, \frak a_{M}^*$, $P({\Bbb A})$, etc, for all 
parabolic subgroups $P$. In particular, one checks that the following 
compactibility condition holds
$$\int_{M_0({\Bbb A})\cdot U_0({\Bbb A})\cdot {\Bbb K}}f(mnk)\,dk\,dn\,
m^{-2\rho_0}dm=\int_{M({\Bbb A})\cdot U({\Bbb A})\cdot {\Bbb K}}f(mnk)\,dk\,
dn\,m^{-2\rho_P}dm$$ for all continuous functions $f$ with compact supports 
on $G({\Bbb A})$, where $\rho_P$ denotes the half of the sum of the positive 
roots of  the maximal split torus $T_P$ of the central $Z_{M}$ of $M$. 
For later use, 
denote also by $\Delta_P$ the set of positive roots determined by $(P,T_P)$ 
and $\Delta_0=\Delta_{P_0}$.

Fix an isomorphism $T_0\simeq {\Bbb G}_m^R$. Embed ${\Bbb R}_+^*$ by the map 
$t\mapsto (1;t)$. Then we obtain a natural injection $({\Bbb R}_+^*)^R
\hookrightarrow T_0({\Bbb A})$ which splits. Denote by $A_{M_0({\Bbb A})}$ 
the unique connected subgroup of $T_0({\Bbb A})$ which projects onto 
$({\Bbb R}_+^*)^R$. More generally, for a standard parabolic subgroup 
$P=MU$, set $A_{M({\Bbb A})}:=
A_{M_0({\Bbb A})}\cap Z_{M({\Bbb A})}$ where as used above $Z_*$ denotes 
the center of the group $*$. Clearly, $M({\Bbb A})=A_{M({\Bbb A})}\cdot 
M({\Bbb A})^1$. For later use, set also 
$A_{M({\Bbb A})}^G:=\{a\in A_{M({\Bbb A})}:\log_Ga=0\}.$ 
Then $A_{M({\Bbb A})}=A_{G({\Bbb A})}\oplus
A_{M({\Bbb A})}^G.$

Note that ${\Bbb K}$, $M(F)\backslash M({\Bbb A})^1$ and $U(F)\backslash 
U({\Bbb A})$ are all compact, thus with the Langlands decomposition 
$G({\Bbb A})=U({\Bbb A})M({\Bbb A}){\Bbb K}$ in mind, the reduction theory 
for $G(F)\backslash G({\Bbb A})$ or more generally 
$P(F)\backslash G({\Bbb A})$ is reduced to that for 
$A_{M({\Bbb A})}$ since $Z_G(F)\cap Z_{G({\Bbb A})}\backslash 
Z_{G({\Bbb A})}\cap G({\Bbb A})^1$ is compact as well. As such 
for $t_0\in 
M_0({\Bbb A})$ set
$$A_{M_0({\Bbb A})}(t_0):=\{a\in A_{M_0({\Bbb A})}:a^\alpha>t_0^\alpha
\forall\alpha\in\Delta_0\}.$$ Then, for a fixed compact subset 
$\omega\subset P_0({\Bbb A})$, we have the corresponding Siegel set 
$$S(\omega;t_0):=\{p\cdot a\cdot k:p\in \omega, a\in A_{M_0({\Bbb A})}
(t_0),k\in {\Bbb K}\}.$$ In particular, for big enough $\omega$ and small 
enough $t_0$, i.e, $t_0^\alpha$ is very close to 0 for all 
$\alpha\in\Delta_0$, the classical reduction theory may be restated as 
$G({\Bbb A})=G(F)\cdot S(\omega;t_0)$. More generally set 
$$A_{M_0({\Bbb A})}^P(t_0):=\{a\in A_{M_0({\Bbb A})}:a^\alpha>t_0^\alpha
\forall\alpha\in\Delta_0^P\},$$ and $$S^P(\omega;t_0):=\{p\cdot a\cdot k:
p\in \omega, a\in A_{M_0({\Bbb A})}^P(t_0),k\in {\Bbb K}\}.$$  
Then similarly as above
for big enough $\omega$ and small enough $t_0$, $G({\Bbb A})=P(F)\cdot 
S^P(\omega;t_0)$. (Here 
$\Delta_0^P$ denotes the set of positive roots for $(P_0\cap M,T_0)$.)

Fix an embedding $i_G:G\hookrightarrow SL_n$ sending $g$ to $(g_{ij})$. 
Introducing a hight function on $G({\Bbb A})$ by setting 
$\|g\|:=\prod_{v\in S}\text{sup}\{|g_{ij}|_v:\forall i,j\}$. It is well-known 
that up to $O(1)$, hight functions are unique. This implies that the following 
growth conditions do not depend on the height function we choose.

A function $f:G({\Bbb A})\to {\Bbb C}$ is said to have moderate growth if 
 there exist 
$c,r\in {\Bbb R}$ such that $|f(g)|\leq c\cdot \|g\|^r$ for all
 $g\in G({\Bbb A})$. Similarly, for a standarde parabolic subgroup $P=MU$, 
a function $f:U({\Bbb A})M(F)\backslash G({\Bbb A})\to {\Bbb C}$ is said 
to have moderate growth if there exist $c,r\in {\Bbb R},\lambda\in 
\text{Re}X_{M_0}$ such that for any $a\in A_{M({\Bbb A})},k\in {\Bbb K}, 
m\in M({\Bbb A})^1\cap S^P(\omega;t_0)$,
$$|f(amk)|\leq c\cdot \|a\|^r\cdot m_{P_0}(m)^\lambda.$$

Also a function $f:G({\Bbb A})\to {\Bbb C}$ is said to be smooth if for any 
$g=g_f\cdot g_\infty \in G({\Bbb A}_f)\times G({\Bbb A}_\infty)$, there exist 
open neighborhoods $V_*$ of $g_*$ in $G({\Bbb A})$ and a $C^\infty$-function 
$f':V_\infty\to {\Bbb C}$ such that $f(g_f'\cdot g_\infty')=f'(g_\infty')$ for 
all $g_f'\in V_f$ and $g_\infty'\in V_\infty$.

By contrast, a function $f: S(\omega;t_0)\to {\Bbb C}$ is said to be rapidly 
decreasing if there exists $r>0$ and for all $\lambda\in \text{Re}X_{M_0}$ 
there exists $c>0$ such that for $a\in A_{M({\Bbb A})}, g\in G({\Bbb A})^1\cap 
S(\omega;t_0)$, $|\phi(ag)|\leq c\cdot\|a\|\cdot m_{P_0}(g)^\lambda$. And a 
function $f:G(F)\backslash G({\Bbb A})\to {\Bbb C}$ is said to be rapidly 
decreasing if $f|_{S(\omega;t_0)}$ is so.

By definition,  a function 
$\phi:U({\Bbb A})M(F)\backslash G({\Bbb A})\to {\Bbb C}$ is called automorphic 
if

\noindent
(i) $\phi$ has moderate growth;

\noindent
(ii) $\phi$ is smooth;

\noindent
(iii) $\phi$ is ${\Bbb K}$-finite, i.e, the ${\Bbb C}$-span of all 
$\phi(k_1\cdot *\cdot k_2)$ parametrized by $(k_1,k_2)\in {\Bbb K}\times 
{\Bbb K}$ is finite dimensional; and 

\noindent
(iv) $\phi$ is $\frak z$-finite, i.e, the ${\Bbb C}$-span of all 
$\delta(X)\phi$ parametrized by all $X\in \frak z$ is finite dimensional. Here 
$\frak z$ denotes the center of the universal enveloping algebra 
$\frak u:=\frak U(\text{Lie}G({\Bbb A}_\infty))$ of the Lie algebra of  
$G({\Bbb A}_\infty)$ and $\delta(X)$
denotes the derivative of $\phi$ along $X$.

For such a function $\phi$, set $\phi_k:M(F)\backslash M({\Bbb A})\to 
{\Bbb C}$ by $m\mapsto m^{-\rho_P}\phi(mk)$ for all $k\in {\Bbb K}$. 
Then one checks that $\phi_k$ is an automorphic form in the usual sense. Set
 $A(U({\Bbb A})M(F)\backslash G({\Bbb A}))$ be the space of automorphic 
forms on $U({\Bbb A})M(F)\backslash G({\Bbb A})$.

For a measurable locally $L^1$-function $f:U(F)\backslash G({\Bbb A})\to 
{\Bbb C}$ define its constant term along with the standard parabolic subgroup 
$P=UM$ to be the function $f_P:U({\Bbb A})\backslash G({\Bbb A})\to {\Bbb C}$ 
given by
$g\to\int_{U(F)\backslash G({\Bbb A})}f(ng)dn.$ 
Then an automorphic form $\phi\in A(U({\Bbb A})M(F)\backslash G({\Bbb A}))$ is 
called a cusp form if for any standard parabolci subgroup $P'$ properly 
contained in $P$, $\phi_{P'}\equiv 0$. Denote by
 $A_0(U({\Bbb A})M(F)\backslash G({\Bbb A}))$  the space of cusp 
forms on $U({\Bbb A})M(F)\backslash G({\Bbb A})$. One checks easily that

\noindent
(i) all cusp forms are rapidly decreasing; and hence

\noindent
(ii) there is a natural pairing
 $$\langle\cdot,\cdot\rangle:A_0(U({\Bbb A})M(F)\backslash G({\Bbb A}))\times 
A(U({\Bbb A})M(F)\backslash G({\Bbb A}))\to {\Bbb C}$$ defined by
$\langle \psi,\phi\rangle:=\int_{Z_{M({\Bbb A})}U({\Bbb A})M(F)\backslash 
G({\Bbb A})}\psi(g)\bar\phi(g)\,dg.$

Moreover, for a (complex) character $\xi:Z_{M({\Bbb A})}\to {\Bbb C}^*$ of 
$Z_{M({\Bbb A})}$ set 
$$\eqalign{A&(U({\Bbb A})M(F)\backslash G({\Bbb A}))_\xi\cr
:=&\{\phi\in A(U({\Bbb A})M(F)
\backslash G({\Bbb A})):\phi(zg)=z^{\rho_P}\cdot\xi(z)\cdot\phi(g),\forall 
z\in Z_{M({\Bbb A})}, g\in G({\Bbb A})\}\cr}$$ and 
$$A_0(U({\Bbb A})M(F)\backslash G({\Bbb A}))_\xi:=A_0(U({\Bbb A})M(F)
\backslash G({\Bbb A}))\cap A(U({\Bbb A})M(F)\backslash G({\Bbb A}))_\xi.$$

Set now $$A(U({\Bbb A})M(F)\backslash G({\Bbb A}))_Z:=\sum_{\xi\in 
\text{Hom}(Z_{M({\Bbb A})},{\Bbb C}^*)}A(U({\Bbb A})M(F)\backslash 
G({\Bbb A}))_\xi$$ and
$$A_0(U({\Bbb A})M(F)\backslash G({\Bbb A}))_Z:=\sum_{\xi\in 
\text{Hom}(Z_{M({\Bbb A})},{\Bbb C}^*)}A_0(U({\Bbb A})M(F)\backslash 
G({\Bbb A}))_\xi.$$ One checks that the natural morphism $${\Bbb C}
[\text{Re}\frak a_M]\otimes A(U({\Bbb A})M(F)\backslash G({\Bbb A}))_Z
\to A(U({\Bbb A})M(F)\backslash G({\Bbb A}))$$ defined by 
$(Q,\phi)\mapsto \big(g\mapsto Q(\log_M(m_P(g))\big)\cdot \phi(g)$ is an 
isomorphism, using the special structure of $A_{M({\Bbb A})}$-finite 
functions and the Fourier analysis over the compact space $A_{M({\Bbb A})}
\backslash Z_{M({\Bbb A})}$. Consequently, we also obtain a natural isomorphism
$${\Bbb C}[\text{Re}\frak a_M]\otimes A_0(U({\Bbb A})M(F)\backslash 
G({\Bbb A}))_Z\to A_0(U({\Bbb A})M(F)\backslash G({\Bbb A}))_\xi.$$

Set also $\Pi_0(M({\Bbb A}))_\xi$ be isomorphism classes of irreducible 
representations of $M({\Bbb A})$ occuring in the space $A_0(M(F)\backslash 
M({\Bbb A}))_\xi$, and
$$\Pi_0(M({\Bbb A}):=\cup_{\xi\in \text{Hom}(Z_{M({\Bbb A})},{\Bbb C}^*)}
\Pi_0(M({\Bbb A}))_\xi.$$ (More precisely, we should 
use $M({\Bbb A}_f)\times (M({\Bbb A})\cap {\Bbb K},
\text{Lie}(M({\Bbb A}_\infty))\otimes_{\Bbb R}{\Bbb C}))$ instead of $M({\Bbb A})$.) For any 
$\pi\in \Pi_0(M({\Bbb A}))_\xi$ set $A_0(M(F)\backslash M({\Bbb A})_\pi$ to 
be the isotypic component of type $\pi$ of $A_0(M(F)\backslash 
M({\Bbb A})_\xi$, i.e, the set of cusp forms of $M({\Bbb A})$ generating a 
semi-simple isotypic $M({\Bbb A}_f)\times (M({\Bbb A})\cap {\Bbb K},
\text{Lie}(M({\Bbb A}_\infty))\otimes_{\Bbb R}{\Bbb C}))$-module of type $\pi$.
Set $$\eqalign{A_0&(U({\Bbb A})M(F)\backslash G({\Bbb A}))_\pi\cr
:=&\{\phi\in 
A_0(U({\Bbb A})M(F)\backslash G({\Bbb A})):\phi_k\in A_0(M(F)\backslash 
M({\Bbb A}))_\pi,\forall k\in {\Bbb K}\}.\cr}$$ Clearly
$$A_0(U({\Bbb A})M(F)\backslash G({\Bbb A}))_\xi=\oplus_{\pi\in 
\Pi_0(M({\Bbb A}))_\xi} A_0(U({\Bbb A})M(F)\backslash G({\Bbb A}))_\pi.$$

More generally, let $V\subset A(M(F)\backslash M({\Bbb A}))$ be an irreducible 
$M({\Bbb A}_f)\times (M({\Bbb A})\cap {\Bbb K},\text{Lie}(M({\Bbb A}_\infty))
\otimes_{\Bbb R}{\Bbb C}))$-module with $\pi_0$ the induced representation of 
$M({\Bbb A}_f)\times (M({\Bbb A})\cap {\Bbb K},\text{Lie}(M({\Bbb A}_\infty))
\otimes_{\Bbb R}{\Bbb C}))$. Then we call $\pi_0$ an automorphic 
representation of $M({\Bbb A})$. Denote by $A(M(F)\backslash 
M({\Bbb A})_{\pi_0}$  the isotypic 
subquotient module of type $\pi_0$ of $A(M(F)\backslash M({\Bbb A})$. One checks
that
$$V\otimes \text{Hom}_{M({\Bbb A}_f)\times (M({\Bbb A})\cap {\Bbb K},
\text{Lie}(M({\Bbb A}_\infty))\otimes_{\Bbb R}{\Bbb C}))}(V,A(M(F)
\backslash M({\Bbb A})))\simeq A(M(F)\backslash M({\Bbb A}))_{\pi_0}.$$ 
Set $$\eqalign{A&(U({\Bbb A})M(F)\backslash G({\Bbb A}))_{\pi_0}\cr
:=&\{\phi\in 
A(U({\Bbb A})M(F)\backslash G({\Bbb A})):\phi_k\in A(M(F)\backslash 
M({\Bbb A}))_{\pi_0},\forall k\in {\Bbb K}\}.\cr}$$
Moreover if 
$A(M(F)\backslash M({\Bbb A}))_{\pi_0}\subset
A_0(M(F)\backslash M({\Bbb A}))$, we call $\pi_0$ a cuspidal representation.

Two automorphic representations $\pi$ and $\pi_0$ of $M({\Bbb A})$ are said 
to be equivalent if there exists $\lambda\in X_M^G$ such that $\pi\simeq 
\pi_0\otimes\lambda$. This, in practice, means that $A(M(F)\backslash 
M({\Bbb A}))_\pi
=\lambda\cdot A(M(F)\backslash M({\Bbb A}))_{\pi_0}.$ That is for any 
$\phi_\pi\in A(M(F)\backslash M({\Bbb A}))_\pi$ there exists a 
$\phi_{\pi_0}\in A(M(F)\backslash M({\Bbb A}))_{\pi_0}$ such that 
$\phi_\pi(m)=m^\lambda\cdot \phi_{\pi_0}(m)$. Consequently, 
$$A(U({\Bbb A})M(F)\backslash G({\Bbb A}))_\pi=(\lambda\circ m_P)\cdot 
A(U({\Bbb A})M(F)\backslash G({\Bbb A}))_{\pi_0}.$$ Denote by 
$\frak P:=[\pi_0]$ the equivalence class of $\pi_0$. Then $\frak P$ is an 
$X_M^G$-principal homogeneous space, hence admits a natural complex structure. 
Usually we call $(M,\frak P)$ a cuspidal datum of $G$ if $\pi_0$ is cuspidal. 
Also for $\pi\in\frak P$ set $\text{Re}\pi:=\text{Re}\chi_\pi=|\chi_\pi|\in 
\text{Re}X_M$, where $\chi_\pi$ is the central character of $\pi$, and 
$\text{Im}\pi:=\pi\otimes(-\text{Re}\pi)$.

Now fix an irreducible automorphic representation $\pi$ of $M({\Bbb A})$ and
$\phi\in A(U({\Bbb A})M(F)\backslash G({\Bbb A}))_\pi$, define the 
associated Eisenstein series $E(\phi,\pi):G(F)\backslash G({\Bbb A})\to 
{\Bbb C}$ by $$E(\phi,\pi)(g):=\sum_{\delta\in P(F)\backslash G(F)}\phi
(\delta g).$$ Then one checks that there is an open cone ${\Cal C}\subset 
\text{Re}X_M^G$ such that if $\text{Re}\pi\in {\Cal C}$, $E(\lambda\cdot \phi,
\pi\otimes\lambda)(g)$ converges uniformly for $g$ in a compact subset of 
$G({\Bbb A})$ and $\lambda$ in an open neighborhood of 0 in $X_M^G$. For 
example, if $\frak P=[\pi]$ is cuspidal, we may even take ${\Cal C}$ to be 
the cone $\{\lambda\in \text{Re}X_M^G:\langle\lambda-\rho_P,
\alpha^\vee\rangle>0,\forall\alpha\in \Delta_P^G\}$.
As a direct consequence, then $E(\phi,\pi)\in A(G(F)\backslash G({\Bbb A}))$.
 That is, it
is an automorphic form.

As noticed above, being an automorphic form, $E(\phi,\pi)$ is of moderate growth.
 However, in general it is not integrable over $Z_{G({\Bbb A})}G(F)\backslash 
G({\Bbb A})$. To remedy this, classically, as initiated in the so-called 
Rankin-Selberg method, analytic truncation is used: From Fourier analysis, 
we understand 
that the probelmatic terms are the so-called constant terms, which are of 
slow growth, so by cutting off them, the reminding one is of rapidly growth 
and hence integrable.

In general, it is very difficult to make such an analytic truncation 
intrinsically related with arithmetic properties of number fields. 
(See however, the Rankin-Selberg method [Bu1,2], [Z] and the Arthur-Selberg trace 
formula [Ar1-4].) 
On the other hand, Eisenstein series themselves are quite intrinsic 
arithmetical invariants.
Thus it is natural for us on one hand to keep Eisenstein series unchanged while on the other 
to find new moduli spaces, which themselves are intrinsically parametrized certain 
modular objects, and over which Eisenstein series are integrable.

This is exactly the approach we take. As said, we are going to view Eisenstein
 series as something globally defined, and use a geo-arithmetical truncation 
for the space $G(F)\backslash G({\Bbb A})$ so that the integrations of the 
Eisenstein series over the newly obtained moduli spaces give us naturally 
non-abelian $L$ functions for  number fields.

As such, let us now come back to the group $G=GL_r$, then as in 2.1, we 
obtain the moduli space ${\Cal M}_{F,r}^{\leq p}[|\Delta_F|^{r\over 2}]$ 
and also a well-defined integration $$L_{F,r}^{\leq p}(\phi,\pi):=
\int_{{\Cal M}_{F,r}^{\leq p}[|\Delta_F|^{r\over 2}]}E(\phi,\pi)(g)\,dg,
\qquad \text {Re}\pi\in {\Cal C}.$$
\eject
\vskip 0.30cm
\noindent
{\bf III.2.3. New Non-Abelian $L$-Functions}

However,  in general, we do not know whether the above  
defined integration has any nice properties such as meromorphic continuation 
and functional equations etc... It is to remedy this that 
we make a further choice of automorphic forms.
  
Fix then a convex polygon $p:[0,r]\to {\Bbb R}$ as in 2.1 and hence 
we get the moduli space ${\Cal M}_{F,r}^{\leq p}[|\Delta_F|^{r\over 2}]$. 
For $G=GL_r$, fix the minimal parabolic subgroup $P_0$ corresponding to the 
partition $(1,\cdots,1)$ with $M_0$ consisting of diagonal matrices. Fix a 
standard parabolic subgroup $P_I=U_IM_I$ corresponding to the partition 
$I=(r_1,\cdots,r_{|P|})$ of $r$ with $M_I$ the standard Levi and $U_I$ the 
unipotent radical.
 
Then for a fixed irreducible automorphic representation $\pi$ of 
$M_I({\Bbb A})$, choose $$\eqalign{\phi\in A&(U_I({\Bbb A})M_I(F)\backslash 
G({\Bbb A}))_\pi\cap 
L^2(U_I({\Bbb A})M_I(F)\backslash G({\Bbb A}))\cr
:=&A^2(U_I({\Bbb A})M_I(F)
\backslash G({\Bbb A}))_\pi,\cr}$$ where $L^2(U_I({\Bbb A})M_I(F)\backslash 
G({\Bbb A}))$ denotes the space of $L^2$ functions on the space 
$Z_{G({\Bbb A})}
U_I({\Bbb A})M_I(F)\backslash G({\Bbb A})$. Denote the associated Eisenstein 
series by $E(\phi,\pi)\in A(G(F)\backslash G({\Bbb A}))$.
\vskip 0.30cm
\noindent
{\bf Main Definition B.} {\it With the same notation as above, a rank $r$ 
non-abelian $L$-function $L_{F,r}^{\leq p}(\phi,\pi)$ for the number field $F$
associated to an
$L^2$-automorphic form $\phi\in A^2(U_I({\Bbb A})M_I(F)\backslash 
G({\Bbb A}))_\pi$ is defined by the following integration
$$L_{F,r}^{\leq p}(\phi,\pi):=
\int_{{\Cal M}_{F,r}^{\leq p}[|\Delta_F|^{r\over 2}]}E(\phi,\pi)(g)\,dg,
\qquad \text {Re}\,\pi\in {\Cal C}.$$}

More generally, for any standard parabolic subgroup $P_J=U_JM_J\supset P_I$ 
(so that the partition $J$ is a refinement of $I$), we have the corresponding 
relative Eisenstein series $$E_I^J(\phi,\pi)(g):=\sum_{\delta\in P_I(F)
\backslash P_J(F)}\phi(\delta g),\forall g\in P_J(F)\backslash G({\Bbb A}).$$ 
It is well-known that there is  an open cone ${\Cal C}_I^J$
in $\text {Re}X_{M_I}^{P_J}$ such that for $\text {Re}\pi\in {\Cal C}_I^J$, 
$E_I^J(\phi,\pi)\in A(P_J(F)\backslash G({\Bbb A})).$ Here $X_{M_I}^{P_J}$ 
is defined similarly as $X_M^G$ with $G$ replaced by $P_J$. Then we have 
a well-defined relative non-abelian $L$-function
$$L_{F,r}^{P_J;\leq p}(\phi,\pi):=
\int_{{\Cal M}_{F,r}^{P_J;\leq p}[|\Delta_F|^{r\over 2}]}E_I^J(\phi,\pi)(g)\,dg,
\qquad \text {Re}\pi\in {\Cal C}_I^J.$$ 
 \vskip 0.30cm
\noindent
{\it Remark.} Here when defining non-abelian $L$-functions we assume that 
$\phi$ comes from a single irreducible automorphic representations, but this 
restriction is rather artifical and can be removed easily, since
such a restriction only serves the purpose of giving the 
constructions and results in a very neat way. 

We end this section by pointing out that the discussion for non-abelian 
$L$-functions holds  for the just defined relative
 non-abelian $L$-functions as well. So from now on, we will leave such a 
 discussion to the reader while concentrate ourselves to non-abelian $L$-functions.
\vskip 0.30cm
\noindent
{\bf III.2.4. Meromorphic Extension and Functional Equations}

With the same notation as above, set $\frak P=[\pi]$. For $w\in W$ the Weyl 
group of $G=GL_r$, fix once and for all representative $w\in G(F)$ of $w$. Set 
$M':=wMw^{-1}$ and denote the associated parabolic subgroup by $P'=U'M'$. 
$W$ acts naturally on the automorphic representations, from which we obtain 
an equivalence
classes $w{\frak P}$ of automorphic representations of $M'({\Bbb A})$. As 
usual, define the associated intertwining operator $M(w,\pi)$ by 
$$(M(w,\pi)\phi)(g):=\int_{U'(F)\cap wU(F)w^{-1}\backslash U'({\Bbb A})}
\phi(w^{-1}n'g)dn',\qquad
\forall g\in G({\Bbb A}).$$ One checks that if $\langle \text {Re}\pi,
\alpha^\vee\rangle\gg 0,\forall \alpha\in\Delta_P^G$,

\noindent
(i) for a fixed $\phi$, $M(w,\pi)\phi$ depends only on the double coset 
$M'(F)wM(F)$. So $M(w,\pi)\phi$ is well-defined for $w\in W$;

\noindent
(ii) the above integral converges absolutely and uniformly for $g$ varying 
in a compact subset of $G({\Bbb A})$;

\noindent
(iii) $M(w,\pi)\phi\in A(U'({\Bbb A})M'(F)\backslash G({\Bbb A}))_{w\pi}$;
and if $\phi$ is $L^2$, which from now on we always assume, so is  
$M(w,\pi)\phi$.
\eject
\vskip 0.30cm
\noindent
{\bf Basic Facts of Non-Abelian $L$-Functions.} {\it With the same notation 
above, 

\noindent
(I) (Meromorphic Continuation) $L_{F,r}^{\leq p}(\phi,\pi)$ for 
$\text{Re}\pi\in {\Cal C}$ is well-defined and admits a unique 
meromorphic continuation to the whole space $\frak P$;

\noindent
(II) (Functional Equation) As meromorphic functions on $\frak P$,
$$L_{F,r}^{\leq p}(\phi,\pi)=L_{F,r}^{\leq p}(M(w,\pi)\phi,w\pi),\qquad \forall w\in W.$$}

\noindent
{\it Proof.} This is a direct consequence of the fundamental results of 
Langlands on Eisenstein series and spectrum decompositions. [See e.g, [Ar1], [La1,2] and [MW1]).
Indeed, if $\phi$ is cuspidal, by definition, (I) is a direct 
consequence of Prop. II.15, Thm. IV.1.8 of [MW] and (II) is a direct 
consequence of Thm. IV.1.10 of [MW].

More generally, if $\phi$ is only $L^2$, then by Langlands' theory of 
Eisenstein series and spectral decomposition, $\phi$ may be obtained as 
the residue of
relative Eisenstein series coming from cuspidal forms, since $\phi$ is $L^2$ 
automorphic. As such then (I) and (II) are direct consequences of the proof 
of VI.2.1(i) at p.264 of [MW].
\vskip 0.30cm
\noindent
{\bf III.2.5. Holomorphicity and Singularities}

Let $\pi\in \frak P$ and $\alpha\in\Delta_M^G$. Define the function 
$h:\frak P\to {\Bbb C}$ by $\pi\otimes\lambda\mapsto \langle \lambda,
\alpha^\vee\rangle,\forall\lambda\in X_M^G\simeq \frak a_M^G$. Here as usual, 
$\alpha^\vee$ denotes the coroot associated to $\alpha$. Set 
$H:=\{\pi'\in\frak P:h(\pi')=0\}$ and call it a root hyperplane. Clearly the 
function $h$ is determined by $H$, hence we also denote  $h$ by $h_H$. 
Note also that root hyperplanes depend on the base point $\pi$ we choose.

Let $D$ be a set of root hyperplanes. Then 

\noindent
(i) the singularities of a meromorphic function $f$ on $\frak P$ is said to 
be carried out by $D$ if for all $\pi\in\frak P$, there exist $n_\pi:D\to 
{\Bbb Z}_{\geq 0}$ zero almost everywhere such that $\pi'\mapsto 
\big(\Pi_{H\in D}h_H(\pi')^{n_\pi(H)}\big)\cdot f(\pi')$ is holomorphic at
$\pi'$;

\noindent
(ii) the singularities of $f$ are said to be without multiplicity at $\pi$ if 
$n_\pi\in\{0,1\}$;

\noindent
(iii) $D$ is said to be locally finite,  if for any compact subset 
$C\subset\frak P$, $\{H\in D:H\cap C\not=\emptyset\}$ is finite.
\vskip 0.30cm
\noindent
{\bf Basic Facts of Non-Abelian $L$-Functions.} {\it With the same notation 
above, 

\noindent
(III) (Holomorphicity) (i) When $\text{Re}\pi\in {\Cal C}$, 
$L_{F,r}^{\leq p}(\phi,\pi)$ is holomorphic;

\noindent
(ii) $L_{F,r}^{\leq p}(\phi,\pi)$ is holomorphic at $\pi$ where 
$\text{Re}\pi=0$;

\noindent
(IV) (Singularities) Assume further that $\phi$ is a cusp form. Then

\noindent
(i) There is a locally finite set of root hyperplanes 
$D$ such that the singularities of $L_{F,r}^{\leq p}(\phi,\pi)$ are carried 
out by $D$;

\noindent
(ii) The singularities of $L_{F,r}^{\leq p}(\phi,\pi)$ are without 
multiplicities at $\pi$ if $\langle \text{Re}\pi,\alpha^\vee\rangle\geq 0,
\forall \alpha\in\Delta_M^G$;

\noindent
(iii) There are only finitely many of singular hyperplanes of 
$L_{F,r}^{\leq p}(\phi,\pi)$ which intersect $\{\pi\in\frak P:
\langle\text{Re}\pi,\alpha^\vee\rangle\geq 0,\forall\alpha\in\Delta_M\}$.}
\vskip 0.30cm
\noindent
{\it Proof.} As above, this is a direct consequence of the fundamental results of 
Langlands on Eisenstein series and spectrum decompositions. [See e.g, [Ar1], [La1,2] and [MW1]).
Indeed, if $\phi$ is a cusp form, (III.i) is a direct consequence of 
Lemma IV.1.7 of [MW], while  (III.ii)  and (IV) are direct consequence of
Prop. IV.1.11 of [MW].

In general when $\phi$ is only $L^2$ automorphic, then we have to use the 
theory of Langlands to realize $\phi$ as the residue of relative Eisenstein 
series defined using cusp forms. (See e.g., item (5) at p.198 and the second 
half part of p.232 of [MW].) 

As such, (III) and (IV) are direct consequence of the definition of 
residue datum and the compactibility between residue and Eisenstein series 
as stated
 for example under item (3) at p.263 of [MW]. 
\vskip 0.30cm
\noindent 
{\it Remarks.} (1) Since $G=GL_r$, one can write down the functional equations
concretely, and give a much more refined result about the singularities of
the non-abelian $L$-functions, for instance, with the 
use of [MW2] about the  residue of Eisenstein 
series. We discuss this elsewhere.

\noindent
(2) As always, once a general construction is given, then we are also 
interseted in some special yet important examples. Here the story should be 
the same.
\eject
\vskip 0.45cm
\noindent
\centerline {\li Chapter IV. Abelian Part of Non-Abelian L-Functions}
\vskip 0.30cm
\centerline
{\li IV.1. Modified Analytic Truncation}
\vskip 0.30cm
Let $G=GL_r$ and $P_0=M_0U_0$ be the minimal parabolic subgroup corresponding 
to the partition $(1,\cdots,1)$. Let $P_1=M_1U_1$ be a fixed standard parabolic
 subgroup with $M_1$ the standard Levi and $U_1$ the unipotent radical.

For a number field $F$ with discriminant $\Delta_F$, let $\pi$ be an 
irreducible automorphic representation of $M_1({\Bbb A})$. Denote by
 $A^2(U_1({\Bbb A})M_1(F)\backslash G({\Bbb A}))_\pi$ the space of 
 $L^2$-automorphic forms in the isotypic component
  $A(U_1({\Bbb A})M_1(F)\backslash G({\Bbb A}))_\pi$.

Then for a fixed convex polygon $p:[0,r]\to {\Bbb R}$ and any $L^2$-automorphic
form $\phi\in A^2(U_1({\Bbb A})M_1(F)\backslash G({\Bbb A}))$ we have the 
associated non-abelian $L$-function
$$L_{F,r}^{\leq p}(\phi;\pi):=\int_{{\Cal M}_{F,r}^{\leq p}
[|\Delta_F|^{r\over 2}]}E(\phi,\pi)(g)\cdot d\mu(g),\qquad 
\text {Re}\,\pi\in {\Cal C}$$ where $E(\phi,\pi)$ denotes the Eisenstein 
series associated to $\phi$ and ${\Cal C}\subset X_{M_1}^G$ is a certain 
positive cone in the previous chapter over which Eisenstein series 
$E(\phi,\pi)$ converges. Also in Ch. III, we show that $L_{F,r}^{\leq p}
(\phi;\pi)$ admits a meromorphic continuation to the whole space 
$\frak P:=[\pi]$, the $X_{M_1}^G$ homogeneous space consisting of automorphic 
representations equivalent to $\pi$ in which a typical element is 
$\pi\otimes\lambda$ with $\lambda\in X_{M_1}^G$.

Motivated by the discussion in II.3, Arthur's analytic truncation [Ar3] and
 Lafforgue's corresponding work on function fields [Laf], we introduce a 
modified  analytic truncation by
$$\Lambda_pE(\phi,\pi)(g):=\sum_P(-1)^{\text{dim}(A_P/Z_G)}\sum_{\delta\in 
P(F)\backslash G(F)}E_P(\phi,\pi)\cdot{\bold 1}(\bar p^{\delta g}>_Pp).$$ 
Then similar as what Arthur does in his fundamental papers [Ar2,3], see also 
Lafforgue [Laf], where function fields are treated, we have the following
\vskip 0.30cm
\noindent
{\bf Basic Facts about the Modified Analytic Truncation.}

\noindent
(i) $\Lambda_p$ is self-adjoint;

\noindent
(ii) $\Lambda_pE(\phi,\pi)$ is rapidly decreasing.

As a direct consequence, we then conclude that
 the integration $$\int_{Z_{G({\Bbb A})}F(F)\backslash G({\Bbb A})}
 \Lambda_p E(\phi,\pi)\cdot d\mu(g),\qquad\text{Re}\pi\in {\Cal C}$$ is 
 well-defined and admits a unique meromorphic  continuation to the 
 whole space $\frak P$.

Motivated by  Arthur's work on the inner product of truncated Eisenstein series [Ar3,4], 
we call the integration $$\int_{Z_{G({\Bbb A})}G(F)\backslash G({\Bbb A})}\Lambda_p
 E(\phi,\pi)\cdot d\mu(g)$$ the leading term of our non-abelian $L$-function 
 $L_{F,r}^{\leq p}(\phi;\pi)$ (when $p$ approach to infinity in a suitable sense). In particular, 
 for the reason we see in Ch. II, we call the above integration {\it an
 abelian part} of the non-abelian $L$-function $L_{F,r}^{\leq p}(\phi;\pi)$
  when $\phi$ is a cusp form,
 and denote it by $L_{F,r}^{\text{ab},\leq p}(\phi;\pi)$.
\vskip 0.30cm
\centerline
{\li IV.2. Relation with Arthur's Analytic Truncation: II}
\vskip 0.30cm
This section is added after Hoffmann at Durham casted his doubts about a 
wishful formula $\Lambda_p\circ\Lambda_p=\Lambda_p$. Clearly, such a formula
for geo-arithmetical truncation is too good to be true. 
As a matter of fact, this is exactly where the key 
difference between our geo-arithmetical truncation and Arthur's analytic 
truncation lies:
For $T\in \frak a_0^+$ a suitable regular element, following 
Arthur, we have the analytic truncation $\Lambda^T\phi$ for any automorphic 
form $\phi$.  Moreover, by Lemma 1.1 of [Ar3], the constant term of 
$\Lambda^T\phi(x)$ along with any standard parabolic subgroup
$P_1$ is zero unless $\varpi(H_0(x)-T)<0$ for all $\varpi\in\hat\Delta_1$. This
is a miracle since, then by definition, 
$$\Lambda^T\circ \Lambda^T=\Lambda^T.$$ 
As a direct consequence,  
$$\int_{Z_{G({\Bbb A})}G(F)\backslash G({\Bbb A})}\Lambda^T
 E(\phi,\pi)\cdot d\mu(g)=\int_{Z_{G({\Bbb A})}G(F)\backslash G({\Bbb A})}
\Lambda^T{\bold 1}\cdot E(\phi,\pi)\cdot d\mu(g).$$ 
Thus if we use Authur's analytic truncation for 
the constant function ${\bold 1}$ to get a truncated (compact) subset in the 
fundamental domain, the integration of Eisenstein series over this 
(analytically) truncated domain is nothing but the leading term 
$$\int_{Z_{G({\Bbb A})}G(F)\backslash G({\Bbb A})}\Lambda^T
 E(\phi,\pi)\cdot d\mu(g).$$ That is to say, no non-abelian part appears for 
such an analytic truncation. 

At this point, the reader may wonder why two truncations 
$\Lambda_p$ and $\Lambda^T$ look very similar superfically, 
yet they yield quite different aspects of the structure --  on one hand, 
Arthur's analytic truncation beautifully singles out very important abelian 
part, on the other, our geo-arithmetical 
truncation creates new rooms for non-abelian contributions. There are deep 
reasons for this: Arthur's analytic 
truncation may be viewed as the \lq second variation\rq$\ $of 
the geo-arithmetical truncation.  Indeed, for any real cocharacter $T$ of 
$M_0$, introduce an associated polygon $p_T:[0,r]\to {\Bbb R}$ such that 

\noindent
(i) it is affine over $[r',r'+1]$ for all 
$r'=0,\dots, r-1$; and 

\noindent
(ii) $\Big(\Delta^2\Big)p_T(r'-1)=(e_{r'}-e_{r'+1})(T)$ where as usual 
$\{e_j-e_{j+1}\}_{j=1}^{r-1}$ denotes the collection of positive roots 
of $SL_r$ and $\Delta f(x):=f(x+1)-f(x)$. 

\noindent
Then by defintion, $${\bold 1}\Big(\Delta^2(p_P^g)>_P
\Delta^2(p_T)\Big)=\tau_P(H_P(g)-T),$$ where $\tau_P(H(g)-T)$ is (one of) 
Arthur's truncation(s). 
\vskip 0.50cm
\centerline
{\li IV.3. A Formula for the Modified Truncated Eisenstein series}
\vskip 0.30cm
Our aim is to give a precise formula for the leading term of 
$L_{F,r}^{\leq p}(\phi;\pi)$, using an advanced version of Rankin-Selberg 
method, motivated by
the result of Langlands-Arthur on the inner product
of truncated Eisenstein series [Ar3,4]. For this purpose we first give a 
formula 
for our analytic truncation of Eisenstein series. This then leads to the 
consideration of constant terms of
Eisenstein series.

According to the theory of Eisenstein series, see e.g, [Ar1], [La1,2] and 
[MW2] 
for any automorphic form $\phi$, 
it is difficult to find a neat formula for the so-called constant terms of the
 corresponding Eisenstein series $E(\phi,\pi)$ along with arbitrary standard 
 parabolic subgroups. However if $\phi$ is cuspidal, then using Bruhat 
 decomposition,
such constant terms may be easily obtained. So from now on, let us assume that
the automorphic form $\phi$ is indeed cuspidal unless we state clearly the 
otherwise.

So let $\phi\in A_0(U_1({\Bbb A})M_1(F)\backslash G({\Bbb A}))$. For any
 standard parabolic subgroup $P=MU$, by the Bruhat decomposition, we have 
$G(F)=\cup_{w\in W_M\backslash W/W_{M_1}}P_1(F)wP(F)$, a disjoint union. Here
$W_*$ is the Weyl group of $*$.
Moreover, one checks that a system of representatives of $W_M\backslash
W/W_{M_1}$ may be given by $$W_{M_1,M}:=\{w\in W:w(\alpha)>0,\forall\alpha\in
R^+(T_0,M),w^{-1}(\beta)>0,\forall \beta\in R^+(T_0,M')\}$$ where $R^+(T_0,*)$
denotes the set of positive roots related to $(T_0,*)$. As a direct 
consequence,
$$\{w\in W_{M_1,M}:wM_1w^{-1}\subset M\}=W(M_1,M)$$ consisting of element 
$w\in W$ such that $wM_1w^{-1}$ is a standard Levi of $M$ and $w^{-1}(\beta)>0$ 
for all $\beta\in R^+(T_0,M)$.

Therefore, by definition, $$\eqalign{E_P(\phi,\pi)
=&\int_{U(F)\backslash U({\Bbb A})}
\sum_{\delta\in P_1(F)\backslash G(F)}\phi(\delta ng)\,dn\cr
=&\sum_{w\in W_{M_1,M}}\sum_{m\in M(F)\cap wP_1(F)w^{-1}\backslash M(F)}
\int_{U(F)\backslash U({\Bbb A})}\sum_{u\in U(F)\cap m^{-1}wP_1(F)w^{-1}m
\backslash U(F)}\phi(w^{-1}mung)dn.\cr}$$ Note that for all $w\in W_{M_1,M}$ 
and $m\in M(F)$,
we have $$\eqalign{&\int_{U(F)\backslash U({\Bbb A})}\sum_{u\in U(F)\cap m^{-1}wP_1(F)
w^{-1}m\backslash U(F)}\phi(w^{-1}mung)\,dn\cr
=&\int_{U(F)\cap wP_1(F)w^{-1}
\backslash U({\Bbb A})}\phi(w^{-1}nmg)\,dn\cr}$$ which is 0 if $w^{-1}Uw\cap M_1 
\not=\{1\}$, or better, if $wM_1w^{-1}$ is not contained in $M$, e.g, $w\not\in
 W(M_1,M)$
 by the cuspidality of $\phi$. 
 
 On the other hand, if $w\in W(M_1,M)$, then
 $M\cap wP_1w^{-1}$ is the standard parabolic subgroup of $M$ with the standard
  Levi $wM_1w^{-1}$ and $U(F)\cap wP_1(F)w^{-1}=U(F)\cap wU_1(F)w^{-1}$. 
  Therefore,
 we obtain a general formula for the constant term of $E(\phi,\pi)$ along $P$ 
 in terms of $\phi$ and its twists by intertwining operators as follows:
 $$E_P(\phi,\pi)(g)=\sum_{w\in W(M_1,M)}\sum_{m\in M(F)\cap wP_1(F) w^{-1}
 \backslash M(F)}\big(M(w,\pi)\phi(\pi)\big)(mg).$$
  
Therefore, 
$$\eqalign{&\Lambda_pE(\phi,\pi)\cr
=&\sum_{P}(-1)^{\text{dim}A_P/Z_G}
\sum_{\delta\in P(F)\backslash G(F)}E_P(\phi,\pi)(\delta g)\cdot 
{\bold 1}(\bar p^{\delta g}>_Pp)\cr
=&\sum_P(-1)^{\text{dim}A_P/Z_G}\sum_{\delta\in P(F)\backslash G(F)}\cr
&\qquad\sum_{w\in W(M_1,M)}\sum_{\xi\in M(F)\cap wP_1(F)w^{-1}\backslash M(F)}
\big(M(w,\pi)\phi\big)(\xi\delta g)\cdot {\bold 1}(\bar p^{\delta g}>_Pp).\cr}$$
Now for any standard parabolic subgroup $P_2$, set $W(\frak a_1,\frak a_2)$ 
to be the set of distinct isomorphisms from $\frak a_1$ onto
$\frak a_2$ obtained by restricting elements in $W$ to $\frak a_1$, where 
$\frak a_i$ 
denotes $\frak a_{P_i}, i=1,2$
Then one checks by definition 
easily that $W(M_1,M)$ is a union over all $P_2$ of elements $w\in W(\frak a_1,\frak a_2)$ 
such that 

\noindent
(i) $w\frak a_1=\frak a_2\supset \frak a_P$; and

\noindent
(ii) $w^{-1}(\alpha)>0, \forall \alpha\in\Delta_{P_2}^P$.

\noindent
Hence, $$\eqalign{&\Lambda_pE(\phi,\pi)\cr
=&\sum_{P_2}\sum_{w\in W(\frak a_1,\frak a_2),P\supset P_2,w^{-1}(\alpha)>0,\forall \alpha\in \Delta_{P_2}^P}
(-1)^{\text{dim}A_P/Z_G}\cr
&\qquad\sum_{\delta\in P(F)\backslash G(F)} {\bold 1}
(\bar p^{\delta g}>_Pp)\cdot
\sum_{\xi\in M(F)\cap wP_1(F)w^{-1}\backslash M(F)}
\big(M(w,\pi)\phi\big)(\xi\delta g)\cr
=&\sum_{P_2}\sum_{w\in W(\frak a_1,\frak a_2)}(-1)^{\text{dim}A_{P_w}/Z_G}
\sum_{\{P:P_2\subset P \subset P_w,w^{-1}(\alpha)>0,\forall \alpha\in \Delta_{P_2}^P\}}
(-1)^{\text{dim}A_P/A_{P_w}}\cr
&\qquad\sum_{\delta\in P(F)\backslash G(F)} {\bold 1}
(\bar p^{\delta g}>_Pp)\cdot
\sum_{\xi\in M(F)\cap wP_1(F)w^{-1}\backslash M(F)}
\big(M(w,\pi)\phi\big)(\xi\delta g).\cr}$$
where for a given $w$, we define $P_w\supset P$ by the conidition that
$$\Delta_{P_2}^{P_w}=\{\alpha\in\Delta_{P_2}:(w\pi)(\alpha^\vee)>0\}.$$
Therefore, since  $${\bold 1}
(\bar p^{\xi\delta g}>_Pp)={\bold 1}
(\bar p^{\delta g}>_Pp), \forall \delta\in P(F)\backslash G(F),
\xi\in P_2(F)\backslash P(F),$$ we have
$$\eqalign{&\Lambda_pE(\phi,\pi)\cr
=&\sum_{P_2}\sum_{w\in W(\frak a_1,\frak a_2)}(-1)^{\text{dim}A_{P_w}/Z_G}
\sum_{\{P:P_2\subset P \subset P_w,w^{-1}(\alpha)>0,\forall \alpha\in 
\Delta_{P_2}^P\}}
(-1)^{\text{dim}A_P/A_{P_w}}\cr
&\qquad\sum_{\delta\in P(F)\backslash G(F)}\sum_{\delta\in P_2(F)\backslash P(F)} \Big({\bold 1}
(\bar p^{\xi\delta g}>_Pp)\cdot
\big(M(w,\pi)\phi\big)(\xi\delta g)\Big)\cr
=&\sum_{P_2}\sum_{w\in W(\frak a_1,\frak a_2)}(-1)^{\text{dim}A_{P_w}/Z_G}
\sum_{\{P:P_2\subset P \subset P_w,w^{-1}(\alpha)>0,\forall \alpha\in \Delta_{P_2}^P\}}
(-1)^{\text{dim}A_P/A_{P_w}}\cr
&\qquad\sum_{\delta\in P_2(F)\backslash G(F)} \Big({\bold 1}
(\bar p^{\delta g}>_Pp)\cdot
\big(M(w,\pi)\phi\big)(\delta g)\Big)\cr
=&\sum_{P_2}\sum_{\delta\in P_2(F)\backslash G(F)}
\sum_{w\in W(\frak a_1,\frak a_2)}(-1)^{\text{dim}A_{P_w}/Z_G}\big(M(w,\pi)\phi\big)(\xi\delta g)\cr
&\qquad\sum_{\{P:P_2\subset P \subset P_w,w^{-1}(\alpha)>0,\forall \alpha\in \Delta_{P_2}^P\}}
(-1)^{\text{dim}A_P/A_{P_w}}
 {\bold 1}
(\bar p^{\delta g}>_Pp).\cr}$$
Set now $${\bold 1}(P_2;p;w):=
\sum_{\{P:P_2\subset P \subset P_w,w^{-1}(\alpha)>0,\forall \alpha\in \Delta_{P_2}^P\}}
(-1)^{\text{dim}A_P/Z_G}
 {\bold 1}
(\bar p^{\delta g}>_Pp).$$
With this, we obtain the following 
\vskip 0.30cm
\noindent
{\bf Lemma.} {\it With the same notation as above,
$$\Lambda_pE(\phi,\pi)(g)=\sum_{P=MU}\sum_{\delta\in P(F)\backslash G(F)}
\sum_{w\in W(M_1), wM_1w^{-1}=M}\big(M(w,\pi)\phi\big)(\delta g)\cdot 
{\bold 1}(P;p;\pi)(\delta g).$$}
\vskip 0.5cm
\centerline
{\li IV.4. Abelian Part of Non-Abelian L-Functions}
\vskip 0.30cm
By the lemma in the previous section, we have
$$\eqalign{&L_{F,r}^{ab;\leq p}(\phi,\pi)\cr
=&\int_{Z_{G({\Bbb A})}G(F)\backslash
 G({\Bbb A})}\sum_{P=MU}\sum_{\delta\in P(F)\backslash G(F)}
\sum_{w\in W(M_1), wM_1w^{-1}=M}\big(M(w,\pi)\phi\big)(\delta g)\cdot 
{\bold 1}(P;p;\pi)(\delta g)\,dg.\cr}$$ 
From an un-folding trick, it is simply
$$\eqalign{&\sum_{P}\sum_{w\in W(M_1), wM_1w^{-1}=M}
\int_{Z_{G({\Bbb A})}P(F)\backslash G({\Bbb A})}
\Big({\bold 1}(P_2;p;w)(g)\cdot
\big(M(w,\pi)\phi\big)(g)\Big)\,dg\cr
=&\sum_{P}\sum_{w\in W(M_1), wM_1w^{-1}=M}
\int_{Z_{G({\Bbb A})}U({\Bbb A})M(F)\backslash G({\Bbb A})}
\Big({\bold 1}_{P}(P;p;w)(g)\cdot
\big(M(w,\pi)\phi\big)(g)\Big)\,dg\cr}$$ where
as usual
${\bold 1}_{P}(P;p;w)(g):=\int_{U(F)\backslash U({\Bbb A})}{\bold 1}
(P;p;w)(n g)\,dn$ denotes the constant term of ${\bold 1}(P;p;w)(g)$ along $P$.

To evaluate this latest integral, we decomposite it into a double
integrations over $$\eqalign{&\big(Z_{G({\Bbb A})}(Z_{M}(F)\cap Z_{M({\Bbb A})}
\backslash Z_{G({\Bbb A})}\cdot Z_{M({\Bbb A})}^1\big)\times
\big(Z_{G({\Bbb A})} Z_{M({\Bbb A})}^1U({\Bbb A)}M(F)\backslash 
G({\Bbb A})\big)\cr
=&\big(Z_{G({\Bbb A})}\cdot Z_{M({\Bbb A})}^1\backslash 
Z_{M({\Bbb A})}\big)\times
\big(Z_{M({\Bbb A})}U({\Bbb A)}M(F)\backslash G({\Bbb A})\big),\cr}$$ where
$Z_{M({\Bbb A})}^1=Z_{M({\Bbb A})}\cap M({\Bbb A})^1$.
That is to say,
$$\eqalign{L_{F,r}^{ab;\leq p}(\phi,\pi)
=&\sum_{P=MU}\sum_{w\in W(M_1), wM_1w^{-1}=M}
\int_{Z_{M({\Bbb A})}U({\Bbb A)}M(F)
\backslash G({\Bbb A})}dg\cr
&\cdot\int_{Z_{G({\Bbb A})}\cdot Z_{M({\Bbb A})}^1
\backslash Z_{M({\Bbb A})}} \Big({\bold 1}_{P}(P;p;w)(zg)\cdot
\big(M(w,\pi)\phi\big)(zg)\Big)dz.\cr}$$

Note now that since $X_{M_1}^G$ has no torsion, there exists a unique element 
$\pi_0$ of $\frak P:=[\pi]$ whose restriction to $A_{M_1({\Bbb A})}^G$ is 
trivial. This then allows to canonically identified $X_{M_1}^G$ with
  $\frak P$ via $\lambda_\pi\in X_{M_1}^G\mapsto \pi:=\pi_0\otimes
  \lambda_\pi\in \frak P$. Hence without loss of generality, we may simply 
  assume that the  restriction of $\pi$ to $A_{M_1({\Bbb A})}^G$ is trivial. 

Therefore, $$\eqalign{L_{F,r}^{ab;\leq p}(\phi,\pi)
=&\sum_{P=MU}\sum_{w\in W(M_1), 
wM_1w^{-1}=M}\int_{Z_{M({\Bbb A})}U({\Bbb A)}M(F)
\backslash G({\Bbb A})}\big(M(w,\pi)\phi\big)(g)\,dg\cr
&\qquad\cdot\int_{Z_{G({\Bbb A})}
\cdot Z_{M({\Bbb A})}^1\backslash Z_{M({\Bbb A})}} \Big({\bold 1}_{P}(P;p;w)(zg)
\cdot
\Big)z^{\rho_P+w\pi}dz.\cr}$$ However as $g$ may be chosen in $G({\Bbb A})^1$, 
clearly,
the integration $$\int_{Z_{G({\Bbb A})}\cdot Z_{M({\Bbb A})}^1\backslash
 Z_{M({\Bbb A})}} \Big({\bold 1}_{P}(P;p;w)(zg)\cdot
\Big)z^{\rho_P+w\pi}dz$$ is independent of $g$. Denote it by $W(P;p;w;\pi)$. 
As a direct consequence, we obtain the following
\vskip 0.30cm
\noindent
{\bf Basic Fact About the Abelian Part.} {\it  With the same notation as above, for $\phi\in
 A_0(U_1({\Bbb A})M_1(F)\backslash G({\Bbb A}))_\pi$,
$$L_{F,r}^{\text{ab},\leq p}(\phi;\pi)=\sum_{P=MU}\sum_{w\in W(M_1), 
wM_1w^{-1}=M}
\Big(W(P;p;w;\pi)\cdot\langle M(w,\pi)\phi,1\rangle\Big).$$}
\vskip 0.30cm
\noindent
{\it Remark.}  With the simple combinatorial lemma Prop 1.1 of [Ar1], one may 
give a concrete formula for ${\bold 1}(P;p;w)$ and hence its constant term  
${\bold 1}_P(P;p;w)$.
Moreover, with the formula above and the Fourier inversion formula, noticing 
the fact that
the dual of $Z_{G({\Bbb A})}M({\Bbb A})^1\backslash M({\Bbb A})$ is simply $X_M^G$, the structure
of $W(P;p;w;\pi)$  may also be understood. We leave all this to the interested reader.
On the other hand, if $\phi$ is only an $L^2$-automorphic form, such an elegent
formula can hardly be given due to the fact that there is no simple formula 
for the constant terms of the associated Eisenstein series. However, if we 
follow Arthur's method, we could give the \lq leading term\rq$\ $of our
non-abelian $L$-function $L_{F,r}^{\leq p}(\phi,\pi)$, as introduced at the 
begining of this chapter. (We reminder the reader that while our
leading term is well-justified, one may simply it by throwing away 
rapidly decreasing terms as done in Arthur [Ar4].) We will discuss this
elsewhere.
\vskip 0.50cm
\centerline {\li Conclusion Remark}
\vskip 0.30cm
All in all, we feel that we have already accomplished our mission to find 
genuine non-abelian $L$-functions: while the abelian part may be
obtained by using Rankin-Selberg(-Langlands-Arthur) method, for the non-commutative theory, we 
should start with these newly found non-abelian $L$-functions. 
In particular, note that from examples of Chapter II,

\noindent
(1) any single non-abelian $L$-function can be
decomposite as the sum of a rank $r$ part and a lower rank part 
coming from parabolic induction; and

\noindent
(2) the rank $r$ part 
consists of two components, the abelian one and the essential 
non-abelian one.

Thus, to understand the universal structure 
of all $L$-functions, both abelian and non-abelian, we should introduce 
 non-abelian Euler products in such a way that 

\noindent
(i) they quantifies what might be called the non-abelian Reciprocity Law; and 

\noindent
(ii) they give an alternative way to introduce
non-abelian $L$-functions in this paper and hence in rank 1 case 
they degenerate to the standard Euler product.
\vskip 0.30cm
It is for this purpose, motivated by
the non-abelian class field theory for function fields over complex numbers, 
to be established in the Appendix, that we propose to study
algebraic relations among images of the Frobenius classes of absolute Galois 
groups in finite (subquotient) groups in $\Sigma_r$. Here 
$\Sigma_r$ is defined to be the collection of all finite groups $G$ satisfying

\noindent
(1) all irreducible unitary representations of $G$ are of rank $\leq r$;

\noindent
(2) $G$ admits at least one rank $r$ irreducible unitary representation.
\eject
\noindent
{\bf Example.} When $r=1$, $\Sigma_r$ is simply the collection of all
finite abelian groups and hence the algebraic relations we are seeking for are 
simply the commutative law. From our point of view, this is why 

\noindent
(i) for abelian $L$-functions, the standrad Euler product exists; and

\noindent
(ii) for non-abelian $L$-functions of rank $r$, as mentioned above,
 the natural decomposition in terms of parabolic data exists.
\vskip 1.0cm
\noindent
{\bf Acknowledgement.} This work is partially supported by JSPS.
When preparing it in approximately 5 years, I receive many helps from various institutes,
in particular, Kobe Univ, Nagoya Univ and Kyushu Univ. I would like to thank them for providing me
with excellent working conditions. Special thanks also due to Ueno  
for his constant support and encouragement;  to van der Geer for kindly explaining his 
joint work with Schoof when we met at Kyoto Univ which then leads to our own work on the 
geo-arithmetical cohomology presented in Ch.I;
to Deninger and Zagier for their discussion, whose question leads to the examples in 
Ch. II. Finally, this work is completed when 
I visit Nottingham. I would like to thank Fesenko for his help and support.
\vskip 1.5cm
\noindent
\centerline {\li REFERENCES} 
\vskip 0.45cm
\noindent
[Ar1] Arthur, J.  Eisenstein series and the trace formula. {\it Automorphic forms, 
representations and $L$-functions}  pp. 253--274, Proc. Sympos. Pure Math., XXXIII, 
Amer. Math. Soc., Providence, R.I., 1979.
\vskip 0.30cm
\noindent
[Ar2] Arthur, J. A trace formula for reductive groups. I. Terms associated to classes 
in $G({\Bbb Q})$. Duke Math. J. 45  (1978), no. 4, 911--952
\vskip 0.30cm
\noindent
[Ar3] Arthur, J. A trace formula for reductive groups. II. Applications of a truncation 
operator. Compositio Math. 40 (1980), no. 1, 87--121.
\vskip 0.30cm
\noindent
[Ar4] Arthur, J. On the inner product of truncated Eisenstein series. Duke Math. J. 49 (1982), no. 1, 35--70.  
\vskip 0.30cm
\noindent
[BV]  Bombieri, E. \& Vaaler, J. On Siegel's lemma. Invent. Math.
{\bf 73} (1983), no. 1, 11--32.
\vskip 0.30cm
\noindent
[Bo1]  Borel, A. Some finiteness properties of adele groups over
number fields, Publ. Math., IHES, {\bf 16} (1963) 5-30
\vskip 0.30cm
\noindent
[Bo2] Borel, A. {\it Introduction aux groupes arithmetictiques}, Hermann, 1969
\vskip 0.30cm
\noindent 
[Bor] A. Borisov, Convolution structures and arithmetic
cohomology, to appear in Comp. Math.
\vskip 0.30cm
\noindent 
[Bos]  Bost, J.-B. Fibr\'es vectoriels hermitiens, degr\'e
d'Arakelov et polygones canoniques, Appendix A to  Exp. No. {\bf 795}, S\'eminaire
Bourbaki 1994/95, Ast\'erisque  {\bf 237} (1996), 154--161.
\vskip 0.30cm
\noindent
[Bu1] Bump, D. {\it Automorphic forms on $GL(3,{\Bbb R})$, springer LNM {\bf 1083}, 1984}
\vskip 0.30cm
\noindent
[Bu2] Bump, D. The Rankin-Selberg method: a survey. {\it Number theory, trace formulas 
and discrete groups} (Oslo, 1987), 49--109, Academic Press, Boston, MA, 1989.
\vskip 0.30cm
\noindent 
[Co] Connes, A. Trace formula in noncommutative geometry and the zeros of the 
Riemann zeta function.  Sel. math. New ser  {\bf 5}  (1999),  no. 1, 29--106.
\vskip 0.30cm
\noindent 
[De1] Deninger, C.  On the $\Gamma$-factors attached to motives.
Invent. Math. {\bf 104} (1991), no. 2, 245--261.
\vskip 0.30cm
\noindent 
[De2] Deninger, C. Motivic $L$-functions and regularized determinants, in
{\it Proc. Sympos. Pure Math}, {\bf 55}, {\it Motives}, edited by U. Jannsen, S. Kleiman and J.-P.
Serre,  (1994), 707-743
\vskip 0.30cm
\noindent 
[De3] Deninger, C. Some analogies between number theory
and dynamical systems on foliated spaces. {\it Proceedings of the International Congress of
Mathematicians}, Vol. I (Berlin, 1998).  Doc. Math.  1998,  Extra Vol. I, 163--186
\vskip 0.30cm
\noindent 
[Fo] Folland, G.B. {\it A course in abstract harmonic analysis}, Studies in advanced
mathematics, CRC Press, 1995
\vskip 0.30cm
\noindent 
[GS] van der Geer \&  Schoof, R. Effectivity of Arakelov
Divisors and the Theta Divisor of a Number Field, Sel. Math., New ser.
{\bf 6} (2000), 377-398  
\vskip 0.30cm
\noindent 
[Gr1] Grayson, D.R. Reduction theory using semistability. 
Comment. Math. Helv.  {\bf 59}  (1984),  no. 4, 600--634.
\vskip 0.30cm
\noindent 
[Gr2] Grayson, D.R.  Reduction theory using semistability. II.  
Comment. Math. Helv.  {\bf 61}  (1986),  no. 4, 661--676.
\vskip 0.30cm
\noindent
[H] Humphreys, J. {\it Introduction to Lie algebras and representation theory}, GTM {\bf 9}, 1972
\vskip 0.30cm
\noindent
[IT] Imai, K \& Terras, A. The Fourier expansion of Eisenstein series for $\text{GL}(3,{\Bbb Z})$. 
Trans. Amer. Math. Soc. 273 (1982), no. 2, 679--694.
\vskip 0.30cm
\noindent 
[Iw1] Iwasawa, K. Letter to Dieudonn\'e, April 8, 1952, in
 {\it Zeta Functions in Geometry}, edited by N.Kurokawa and T. Sunuda,
Advanced Studies in Pure Math. {\bf 21} (1992), 445-450
\vskip 0.30cm
\noindent 
[Iw2] Iwasawa, K. Lectures notes on Riemann(-Artin) Hypothesis, noted by Kimura, Princeton Unir. 197?
\vskip 0.30cm
\noindent 
[L1]  Lang, S. {\it Algebraic Number Theory}, 
Springer-Verlag, 1986
\vskip 0.30cm
\noindent 
[L2] Lang, S. {\it Fundamentals on Diophantine Geometry}, 
Springer-Verlag, 1983
\vskip 0.30cm
\noindent 
[L3] Lang, S. {\it Introduction to Arekelov Theory}, 
Springer-Verlag, 1988
\vskip 0.30cm
\noindent
[Laf] Lafforgue, L. {\it Chtoucas de Drinfeld et conjecture de Ramanujan-Petersson}. 
Asterisque No. 243 (1997)
\vskip 0.30cm
\noindent
[La1] Langlands, R. Eisenstein series. {\it Algebraic Groups and Discontinuous Subgroups}
 Proc. Sympos. Pure Math., IX, Amer. Math. Soc., Providence, R.I., 1979,  pp. 235--252.
\vskip 0.30cm
\noindent
[La2] Langlands, R. {\it On the functional equations satisfied by Eisenstein series}, 
Springer LNM {\bf 544}, 1976
\vskip 0.30cm
\noindent 
[Li] Li, X.  A note on the Riemann-Roch theorem for
function fields.  {\it Progr. Math.}, 139, (1996),  567--570.
\vskip 0.30cm
\noindent
[MW1] Moeglin, C. \& Waldspurger, J.-L. Le spectre residuel de GL(n), Ann. Sci. Ec. Norm. Sup.
{\bf 22}, (1989), 605-674
\vskip 0.30cm
\noindent
[MW2] Moeglin, C. \& Waldspurger, J.-L. {\it Spectral decomposition and Eisenstein series}. 
Cambridge Tracts in Mathematics, {\bf 113}. Cambridge University Press, Cambridge, 1995.
\vskip 0.30cm
\noindent 
[Mo] Moreno, C. {\it  Algebraic curves over finite fields.}
Cambridge Tracts in Mathematics, {\bf 97}, Cambridge University Press, 1991
\vskip 0.30cm
\noindent 
[Mor] Moriwaki, A. Stable sheaves on arithmetic curves, a personal
note dated in 1992
\vskip 0.30cm
\noindent 
[Mu] Mumford, D. {\it Geometric Invariant Theory}, Springer-Verlag,
(1965)
\vskip 0.30cm
\noindent 
[Neu] Neukirch, {\it Algebraic Number Theory}, Grundlehren der
Math. Wissenschaften, Vol. {\bf 322}, Springer-Verlag, 1999
\vskip 0.30cm
\noindent 
[Pa] Parshin, A.N. On the arithmetic of two-dimensional schemes. I. Distributions
and residues. (Russian)  Izv. Akad. Nauk SSSR Ser. Mat.  {\bf 40}  (1976), no. 4, 736--773, 949.
\vskip 0.30cm
\noindent 
[Se1]  Serre, J.-P.  Zeta and $L$ functions, in {\it Arithmetical Algebraic
Geometry} (Proc. Conf. Purdue Univ., 1963)   Harper \& Row (1965), 82-92
\vskip 0.30cm
\noindent 
[Se2] Serre, J.-P. {\it Algebraic Groups and Class Fields}, GTM
{\bf 117}, Springer-Verlag (1988)
\vskip 0.30cm
\noindent
[SC] Selberg, A. \& Chowla, S. On Epstein's zeta-function, J. Reine und Angew. Math {\bf 227} (1967), 86-110 
\vskip 0.30cm
\noindent 
[S] Siegel, C.L. {\it Lectures on the geometry of numbers}, notes by B.
Friedman, rewritten by K. Chandrasekharan with the assistance of R. Suter, Springer-Verlag,
1989.
\vskip 0.30cm
\noindent 
[St1] Stuhler, U. Eine Bemerkung zur Reduktionstheorie
quadratischer Formen.  Arch. Math.  {\bf 27} (1976), no. 6,
604--610.
\vskip 0.30cm
\noindent 
[St2] Stuhler, U.  Zur Reduktionstheorie der positiven
quadratischen Formen. II.  Arch. Math.   {\bf 28}  (1977), no. 6, 611--619. 
\vskip 0.30cm
\noindent 
[Ta] Tate, J. Fourier analysis in number fields and Hecke's
zeta functions, Thesis, Princeton University, 1950 
\vskip 0.30cm
\noindent
[Te] Terras, A. {\it Harmonic analysis on symmetric spaces and applications II}, Springer-Verlag, 1988
\vskip 0.30cm
\noindent
[V] Venkov, A.B. On the trace formula for $SL(3,{\Bbb Z})$, J. Soviet Math., {\bf 12} (1979), 384-424
\vskip 0.30cm
\noindent
[W1] Weil, A. {\it Adeles and algebraic groups}, Prog. in Math, {\bf 23} (1982)
\vskip 0.30cm
\noindent
[W2] Weil, A. {\it Basic Number Theory}, Springer-Verlag, 1973
\vskip 0.30cm
\noindent
[We1] Weng, L. Non-Abelian Class Field Theory for Riemann Surfaces, at
math.AG/0111240
\vskip 0.30cm
\noindent 
[We2] Weng, L.  A Program for Geometric Arithmetic, at math.AG/0111241
\vskip 0.30cm
\noindent 
[We3] Weng, L. A telefax to Deninger and Fesenko, 2002
\vskip 0.30cm
\noindent
[Z] Zagier, D.  The Rankin-Selberg method for automorphic functions which are not of rapid decay. 
J. Fac. Sci. Univ. Tokyo Sect. IA Math. 28 (1981), no. 3, 415--437 (1982).
\vfill
\eject  
\centerline{\li APPENDIX:}
\vskip 0.45cm
\centerline {\we Non-Abelian Class Field Theory for Function Fields Over C} 
\vskip 0.45cm 
\centerline {\li Lin WENG} 
\vskip 0.30cm 
\centerline {\bf Graduate School of Mathematics, Kyushu University, Japan} 
\vskip 0.50cm
In this appendix, using what we call a micro reciprocity law, we complete 
Weil's program [W] for non-abelian class field theory of Riemann 
surfaces. 
\vskip 0.30cm 
\noindent 
{\li 1. Refined Structures for Tannakian Categories} 
\vskip 0.30cm 
Let ${\Bbb T}$ be a Tannakian category with a fiber functor $\omega:{\Bbb 
T}\to \text{Ver}_{\Bbb C}$, where $\text{Ver}_{\Bbb C}$ denotes the category 
of finite dimensional ${\Bbb C}$-vector spaces. An object $t\in {\Bbb T}$ 
is called {\it reducible} if there exist non-zero objects $x,y\in {\Bbb T}$ 
such that $t=x\oplus y$. An object is called {\it irreducible} if it is 
not reducible. If moreover every object $x$ of ${\Bbb T}$ can be written 
uniquely as a sum of irreducible objects $x=x_1\oplus 
x_2\oplus\dots\oplus x_n$, then ${\Bbb T}$ is called a {\it unique 
factorization} Tannakian category. Usually, we call $x_i$'s the 
{\it irreducible components} of $x$. 

A Tannakian subcategory ${\Bbb S}$ of a unique factorization Tannakian 
category ${\Bbb T}$ is called {\it completed} if for $x\in {\Bbb S}$, 
all its irreducible components $x_i$'s in ${\Bbb T}$ are also in ${\Bbb S}$. 
${\Bbb S}$ is called {\it finitely generated} if as an abelian category, it 
is generated by finitely many objects. Moreover,  ${\Bbb S}$  is called {\it 
finitely completed} if (a) ${\Bbb S}$ is finitely generated; (b) ${\Bbb S}$ is 
completed; and (c) $\text{Aut}^{\otimes}\omega\big|_{\Bbb S}$ is a finite 
group. 
\vskip 0.30cm 
\noindent 
{\li 2. An Example} 
\vskip 0.30cm 
Let $M$ be a compact Riemann surface of genus $g$. Fix an effective 
divisor $D=\sum_{i=1}^Ne_iP_i$ with $e_i\in {\Bbb Z}_{\geq 2}$ once for 
all. For simplicity, in this note we always assume that $(M;D)\not=({\Bbb 
P}^1;e_1P_1)$, or $({\Bbb P}^1;e_1P_1+e_2P_2)$ with $e_1\not=e_2$. (These 
cases may be easily treated.) 

By definition a parabolic semi-stable bundle 

\centerline {$\Sigma:=(E=:E(\Sigma);P_1,\dots,P_N; Fil (E|{P_1}),\dots, 
Fil(E|_{P_N});a_{11},\dots, a_{1r_1};\dots; a_{N1},\dots, a_{Nr_N})$} 

\noindent 
of parabolic degree 0 is called a {\it GA bundle} over $M$ along $D$ if 
(i) the parabolic weights are all rational, i.e., 
$a_{ij}\in {\Bbb Q}\cap [0,1)$; (ii) there exist 
$\alpha_{ij}\in {\Bbb Z}, \beta_{ij}\in {\Bbb Z}_{>0}$ such that (a) 
$(\alpha_{ij},\beta_{ij})=1;$ (b) $a_{ij}=\alpha_{ij}/\beta_{ij};$ and (c) 
$\beta_{ij}|e_i$, for all $i,j$. Denote by 
$[\Sigma]$ the Seshadri equivalence class associated with $\Sigma$. 
Moreover, for $[\Sigma]$, define $\omega_D([\Sigma])$ as $E(\text{
Gr}(\Sigma))|_P$, i.e., the fiber of the bundles associated with 
Jordan-H\"older graded parabolic bundle of $\Sigma$ at a fixed $P\in 
M^0:=M\backslash |D|=M\backslash\{P_1,\dots,P_N\}$. 
\vskip 0.30cm 
\noindent 
{\bf Proposition.} {\it With  the same notation as above, put ${\Cal 
M}(M;D):=\{[\Sigma]:\Sigma\ \text{is\ a\ GA\ bundle\ over}\ M\ \text{along} 
\ D\}$. Then ${\Cal M}(M;D)$ is a unique factorization Tannakian category 
and $\omega_D:{\Cal M}(M;D)\to \text{Vec}_{\Bbb C}$ is a fiber functor.} 
\vskip 0.30cm 
\noindent 
{\it Proof.} (1)  By a result of Mehta-Seshadri [MS, Prop. 1.15], 
${\Cal M}(M;D)$ is an abelian category. Then from the unitary 
representation interpretation of a GA bundle, a fundamental result due to 
Seshadri, (see [MS, Thm 4.1], also in Step 3 of Section 5 below,) 
${\Cal M}(M;D)$ is closed under tensor product. The rigidity may be 
checked directly. 

\noindent 
(2) Since the Jordan-H\"older graded bundle is a direct sum 
of stable and hence irreducible objects and is unique, (see [MS, Rm 
1.16],) so, ${\Cal M}(M;D)$ is  a unique factorization category. 

\noindent 
(3) By definition, we know that the functor $\omega$ is exact and tensor. 
So we should check whether it is faithful. This then is a direct 
consequence of the fact that ${\Cal M}(M;D)$ is  a unique factorization 
category and that any morphism between two irreducible objects is either 
zero or a constant multiple of the identity map. 
\vskip 0.30cm 
\noindent 
{\li 3. Reciprocity Map} 
\vskip 0.30cm 
In $[\Sigma]$, choose its associated Jordan-H\"older graded bundle $\text{
Gr}(\Sigma)$ as a representative. Then by the above mentioned 
fundamental result of Seshadri, 
$\text{Gr}(\Sigma)$ corresponds to a unitary representation $\rho_{\text{
Gr}(\Sigma)}:\pi_1(M^0)\to U(r_\Sigma)$, where $r_\Sigma$ denotes the rank 
of 
$E(\Sigma):=E$. 

For each element $g\in \pi_1(M^0)$, we then obtain a ${\Bbb 
C}$-isomorphism of $E(\text{Gr}(\Sigma))|_P$. Thus, in particular, we get 
a natural morphism $$W:\pi_1(M^0)\to \text{Aut}^\otimes\omega_D.$$ 

Now note that $\pi_1(M^0)$ is generated by $2g$ hyperbolic 
transformations $A_1,B_1,\dots,A_g,B_g$ and $N$ parabolic transformations 
$S_1,\dots,S_N$ satisfying a single relation 
$A_1B_1A_1^{-1}B_1^{-1}\dots A_gB_gA_g^{-1}B_g^{-1}S_1\dots S_N=1,$ and 
that $\rho_{\text{Gr}(\Sigma)}(S_i^{e_i})=1$ for all $i=1,\dots, N$. 
([MS, \S 1].)  Denote by $J(D)$ 
the normal subgroup of $\pi_1(M^0)$ generated by 
$S_1^{e_1},\dots,S_N^{e_N}$. Then naturally we obtain the following 
{\it reciprocity map} 
$$W(D):\pi_1(M^0)/J(D)\to \text{Aut}^\otimes \omega_D.$$ 
\vskip 0.30cm 
\noindent 
{\li 4. Main Theorem} 
\vskip 0.30cm 
As usual, a Galois covering $\pi:M'\to M$ is called branched at most at 
$D$ if (1) $\pi$ is  branched at $P_1,\dots,P_N$; and (2) the 
ramification index $e_i'$ of points over $P_i$ divides $e_i$ for all 
$i=1,\dots, N$. Clearly, by changing $D$, we get all finite Galois 
coverings of $M$. 
\vskip 0.30cm 
\noindent 
{\bf Main Theorem.} (1) (Existence and Conductor Theorem) 
{\it There is a natural one-to-one correspondence $w_D$ between $$\{{\Bbb 
S}:\text{finitely\ completed\ Tannakian\ subcategory\ of} \ {\Cal 
M}(M;D)\}$$ and $$ 
\{\pi:M'\to M:\text{finite\ Galois\ covering\ branched\ at\ most\ at}\ 
D\};$$} 
\noindent 
(2) (Reciprocity Law) {\it There is a natural group isomorphism} 
$$\text{Aut}^\otimes (\omega_D\big|_{\Bbb S})\simeq \text{Gal}\,(w_D({\Bbb 
S})).$$ 
\vskip 0.20cm 
\noindent 
{\li 5. Proof} 
\vskip 0.30cm 
\noindent 
{\it Step 1: Galois Theory.} By a result of Bungaard, Nielsen, Fox, M. 
Kato, and Namba, (see e.g., [Na, Thms 1.2.15 and 1.3.9],) we know 
that the assignment 
$(\pi:M'\to M)\mapsto 
\pi_*(\pi_1(M'\backslash\pi^{-1}\{P_1,\dots,P_N\}))$ gives a one-to-one 
correspondence between isomorphism classes of finite Galois coverings 
$\pi:M'\to M$ branched at most at $D$ and finite index (closed) normal 
subgroups $K=K(\pi)$ of $\pi_1(M^0)$ containing $J(D)$. Moreover, we have 
a natural isomorphism $\text{Gal}(\pi)\simeq \pi_1(M^0)/K(\pi)\Big(\simeq 
\big(\pi_1(M^0)/J(D)\big)\Big/\big(K(\pi)/J(D)\big)\Big).$ 
Thus the problem is transformed to the one for finite index 
normal subgroups of $\pi_1(M^0)$ which contain $J(D)$, or the same, finite 
index normal subgroups of $G(D):=\pi_1(M^0)/J(D)$. 

\noindent 
{\it Step 2: Tannakian Category Theory: Geometric Side.} Consider now the 
category ${\Bbb T}(D)$ 
of equivalence classes of unitary representations of $G(D)$. Clearly, 
${\Bbb T}(D)$ forms a unique factorization Tannakian category, whose fiber 
functor 
$\omega(D)$ may be defined to be the forget functor. Now fixed once for 
all a representative 
$\rho_\Sigma:G(D)\to U(r_\Sigma)$ for each equivalence classes $[\Sigma]$. 
(The choice of the representative will not change the essentials below as 
the resulting groups are isomorphic to each other.) 

Let ${\Bbb S}$ be a finitely completed Tannakian subcategory of ${\Bbb 
T}(D)$. Then as in the definition of reciprocity map above, we have a 
natural morphism $G(D)\buildrel\omega_{\Bbb S}\over \to \text{Aut}^\otimes 
\omega\Big|_{\Bbb S}.$ Denote its kernel by $K({\Bbb S})$. Then, by 
definition, $G(D)/K({\Bbb S})$ is a finite group, and 
$K({\Bbb S})=\cap_{[\Sigma]\in{\Bbb S}}\text{ker}\rho_\Sigma.$ 

Since ${\Bbb S}$ is finitely completed, there exists a finite set 
$S=\{[\Sigma_1],[\Sigma_2],\dots,[\Sigma_t]\}$ which generates ${\Bbb S}$ 
as a completed Tannakian subcategory. Set 
$[\Sigma_0]:=\oplus_{i=1}^t[\Sigma_i]$. Then  for any $[\Sigma]\in {\Bbb 
S}$, $\text{ker}(\rho_{\Sigma_0})\subset \text{ker}(\rho_{[\Sigma]})$, 
since ${\Bbb S}$ is generated by $S$. Also, by definition, 
$\text{ker}(\rho_{\Sigma_0})=\cap_{i=0}^t \text{ker}(\rho_{\Sigma_i})$. 
Thus $K({\Bbb S})=\text{ker}(\rho_{\Sigma_0})$. 

Therefore, for any $[\Sigma]\in {\Bbb S}$, 
$\rho_\Sigma=\tilde\rho_\Sigma\circ\Pi(D;{\Bbb S})$ where $\Pi(D;{\Bbb 
S}):G(D)\to G(D)/K({\Bbb S})$ denotes the natural quotient map and 
$\tilde \rho_\Sigma$ is a suitable unitary representation of $G(D)/K({\Bbb 
S})$. 

Now set $\tilde {\Bbb S}:=\{[\tilde\rho_{\Sigma}]:[\Sigma]\in {\Bbb S}\}$. 
$\tilde {\Bbb S}$ is a finitely completed Tannakian subcategory in 
${\Cal U}ni(G(D)/K({\Bbb S}))$, the category of equivalence classes 
of unitary representations of $G(D)/K({\Bbb S})$. 

In particular, since the unitary 
representation $\tilde\rho_{\Sigma_0}$ of $G(D)/K({\Bbb S})$ maps 
$G(D)/K({\Bbb S})$ injectively into its image,  for any two 
elements $g_1,g_2\in G(D)/K({\Bbb S})$, 
$\tilde\rho_{\Sigma_0}(g_1)\not=\tilde\rho_{\Sigma_0}(g_2)$. 

With this, by applying the van Kampen Completeness Theorem ([Ka]), 
which claims that for any compact group $G$, if Z is a subset of the 
category 
${\Cal U}ni(G)$ such that for any two elements $g_1,g_2$, there exists a 
representation  $\rho_{g_1,g_2}$ in $Z$ such that 
$\rho_{g_1,g_2}(g_1)\not=\rho_{g_1,g_2}(g_2)$, then the completed 
Tannakian subcategory generated by $Z$ is the whole 
category ${\Cal U}ni(G)$ itself, we conclude that $\tilde {\Bbb S}={\Cal 
U}ni(G(D)/K({\Bbb S}))$. But as categories, ${\Bbb S}$ is equivalent to 
$\tilde {\Bbb S}$, thus by the Tannaka duality, (see e.g., [DM] or [Ta]) 
we obtain a natural isomorphism 
$\text{Aut}^\otimes\omega|_{\Bbb S}\simeq G(D)/K({\Bbb S}).$ 

On the other hand, if $K$ is a finite index (closed) normal subgroup of 
$G(D)$. Set $\tilde {\Bbb S}:={\Cal U}ni(G(D)/K)$ with the fiber functor 
$\omega_{\tilde {\Bbb S}}$. Compositing with the natural quotient map 
$\Pi:G(D)\to G(D)/K$ we then obtain an equivalent  category ${\Bbb S}$ 
consisting of corresponding unitary representations of 
$G(D)$. ${\Bbb S}$ may also be viewed as a Tannakian subcategory of ${\Cal 
U}ni(G(D))$. We next show that indeed such an ${\Bbb S}$ is a finitely 
completed Tannakian subcategory. 

From definition, $\text{Aut}^\otimes\omega|_{\Bbb S}\simeq \text{
Aut}^\otimes\omega_{\tilde{\Bbb S}}$ which by the Tannaka duality theorem 
is ismorphic to $G(D)/K$. So it then suffices to show that ${\Bbb S}$ is 
finitely generated. But this is then a direct consequence of the fact that 
for any finite group there always exists a unitary representation such 
that the group is injectively mapped into the unitary group. 

\noindent 
{\it Step 3: A Micro Reciprocity Law.} With Steps 1 and 2, the proof of 
the Main Theorem is then completed by the following 

\noindent 
{\bf Weil-Narasimhan-Seshadri Correspondence} (Seshadri, [MS, Thm 4.1]) 
{\it There is a natural one-to-one 
correspondence between  isomorphism classes of unitary 
representations of fundamental groups of $M^0$ 
and equivalence classes of semi-stable parabolic bundles over $M^0$ of 
parabolic degree zero.} 

Indeed, with this theorem, the Seshadri equivalence classes 
of GA bundles over 
$M$ along 
$D$ correspond  naturally in one-to-one to the equivalence classes of 
unitary representations of the group $\pi_1(M^0)/J(D)$. Thus by Step 2, 
the finitely completed Tannakian subcategories of ${\Cal M}(M;D)$ are in 
one-to-one correspond to the finite index closed normal subgroup of 
$\pi_1(M^0)/J(D)$, which by Step 1 in one-to-one correspond to the finite 
Galois coverings of $M$ branched at most at $D$. This gives the existence 
theorem. Along the same line, we have the reciprocity law as well.
\vskip 0.50cm 
\noindent 
{\li 6. An Application to Inverse Galois Problem} 
\vskip 0.30cm  
As it stands, we may use our main theorem to see whether a finite group can 
be realized as a Galois group of certain coverings by constructing 
suitable finitely completed Tannakian subcategory of ${\Cal 
M}(M;D)$. In particular, to relate with the so-called Inverse Galois 
Problem, or better, the Regular Inverse Galois Problem, via a fundamental 
result of Belyi, we only need to consider the finitely completed 
Tannakian subcategories of ${\Cal 
M}({\Bbb P}^1;\Delta)$, for some effective divisors 
$\Delta=e_\infty\cdot[\infty]+e_1\cdot[1]+e_0\cdot[0]$ supported on 
$\{\infty, 1,0\}\subset {\Bbb P}^1$. 
\vskip 0.30cm 
\noindent 
{\bf An Example.} By definition, a collection of objects of a finitely 
completed Tannakian 
subcategory in ${\Cal M}(M;D)$ is called primitive generators if it 
is a smallest collection which generates the subcategory as an {\it 
abelian category}. 

Since ramification points are all fixed to be $\{\infty, 1,0\}$, to simplify the notation,
set $\Sigma_{11}:=({\Cal O};0;0;0); 
\Sigma_{12}:=({\Cal O}(-1);0;{1\over 2};{1\over 2})$ and $\Sigma_{21}:=({\Cal 
O}(-1)\oplus {\Cal O}(-1);{2\over 3},{1\over 3};{1\over 2},0;{1\over 
2},0)$. Then we may view $\Sigma_{ij}$ as parabolic vector bundles over 
${\Bbb P}^1$ with parabolic points $\infty, 1$ and $0$. That is to say, 
the bundle part is given by the line bundles or rank two vector bundle 
over ${\Bbb P}^1$, while the parabolic weights are given by 
$(0;0;0), (0;{1\over 2};{1\over 2})$ and $({2\over 3},{1\over 3};{1\over 
2},0;{1\over 2},0)$ at points $(\infty;1;0)$. Here for rank two cases, 
the filtration is choosen to be the one such that $E_P\supset {\Bbb 
Q}(e_1+e_2)\supset \{0\}$ with the i-th factor gives ${\Cal O}(-1)_P={\Bbb 
C}\cdot e_i$. 
\vskip 0.30cm 
\noindent 
{\bf Proposition.} {\it Let ${\Cal 
R}(\Sigma_{11},\Sigma_{12},\Sigma_{21})$ be the Tannakian subcategory 
generated by the parabolic 
vector bundles $\Sigma_{11},\Sigma_{12}$ and $\Sigma_{21}$. 
Then ${\Cal 
R}(\Sigma_{11},\Sigma_{12},\Sigma_{21})$ is finitely generated with 
$\{\Sigma_{11},\Sigma_{12}, 
\Sigma_{21}\}$ primitive generators. Moreover,  
$$\text{Aut}^{\otimes}\omega\big|_{{\Cal 
R}(\Sigma_{11},\Sigma_{12},\Sigma_{21})}\simeq S_3.$$} 
\vskip 0.30cm 
\noindent 
{\it Remarks.} (1) The reader should first notice that while ${\Cal R}$ 
is  generated by $\Sigma_{ij}$ as {\it Tannakian subcategory}, ${\Cal R}$ 
as an abelian category is already generated by $\Sigma_{ij}$. 
\vskip 0.30cm 
\noindent 
(2) This proposition, via our non-abelian reciprocity law, reproves a 
well-known result that $S_3$ may be realized as a Galois group 
of a branched covering of ${\Bbb P}^1$ ramified at $\infty, 1,0$ with 
remification index $3,2,2$ respectively. 
\vskip 0.30cm 
\noindent 
{\it Proof.} We first need to show that all tensor products could be 
realized as extensions of $\Sigma_{ij}$. For this, we check case by case as follows. 
\vskip 0.30cm
\noindent 
(1) Clearly, $\Sigma_{11}\otimes \Sigma_{ij}=\Sigma_{ij}$; 

\noindent 
(2) $\Sigma_{12}\otimes \Sigma_{12}=\Sigma_{11}.$
Indeed, $$\eqalign{\Sigma_{12}&\otimes \Sigma_{12}\cr
=&({\Cal O}(-2);0+0;{1\over 2}+{1\over 
2};{1\over 2}+{1\over 2})=({\Cal O}(-2);0;1;1)\cr
=&({\Cal O}(-2+2\cdot 
1);0;1-1;1-1)=({\Cal O};0;0;0)\cr
=&\Sigma_{11};\cr}$$ 

\noindent 
(3) $\Sigma_{12}\otimes \Sigma_{21}=\Sigma_{21}.$
Indeed, $$\eqalign{\Sigma_{12}&\otimes \Sigma_{21}\cr
=&({\Cal O}(-2)\oplus {\Cal O}(-2); 
{2\over 3}+0,{1\over 3}+0;{1\over 2}+{1\over 2},{1\over 2}+0;{1\over 
2}+{1\over 2},{1\over 2}+0)\cr
=&({\Cal O}(-2)\oplus {\Cal O}(-2); 
{2\over 3},{1\over 3};1,{1\over 2};1,{1\over 2})\cr 
=&({\Cal O}(-2+1)\oplus {\Cal O}(-2+1); 
{2\over 3},{1\over 3};{1\over 2};0;{1\over 2};0)\cr
=&\Sigma_{21};\cr}$$ 

\noindent 
(4) $\Sigma_{21}\otimes \Sigma_{21}=\Sigma_{21}\oplus \Sigma_{12}\oplus \Sigma_{11}.$
Indeed, $$\eqalign{\Sigma_{21}&\otimes \Sigma_{21}\cr
=&({\Cal O}(-2)^{\oplus 4};{2\over 
3}+{2\over 3},{2\over 3}+{1\over 3},{1\over 3}+{2\over 3},{1\over 
3}+{1\over 3};{1\over 2}+{1\over 2},{1\over 2}+0,0+{1\over 2},0+0;{1\over 
2}+{1\over 2},{1\over 2}+0,0+{1\over 2},0+0)\cr
=&({\Cal O}(-2)^{\oplus 
4};1+{1\over 3},1,1,{2\over 3};1,{1\over 
2},{1\over 2},0;1,{1\over 2},{1\over 2},0)\cr
=&({\Cal O}(-2+1)^{\oplus 
2}\oplus {\Cal O}(-2+1)\oplus {\Cal O}(-2+2);1+{1\over 
3}-1,1-1,1-1,{2\over 3};1-1,{1\over 2},{1\over 2},0;1-1,{1\over 2},{1\over 
2},0)\cr
=&({\Cal O}(-1)^{\oplus 
2}\oplus {\Cal O}(-1)\oplus {\Cal O}(0);{2\over 
3},{1\over 3},0,0;{1\over 2},{1\over 2},0,0;{1\over 2},{1\over 
2},0,0)\cr
=&\Sigma_{21}\oplus \Sigma_{12}\oplus \Sigma_{11}.\cr}$$ 
 
Now note that the above decompositions for the tensor products are 
exactly the same as that for the 
irreducible representations of $S_3$. So ${\Cal 
R}(\Sigma_{11},\Sigma_{12},\Sigma_{21})$ is equivalent to the category of finite dimensional complex representations of $S_3$.
As a direct consequence, 
we conclude that the structure group $\text{Aut}^{\otimes}\omega\big|_{{\Cal 
R}(\Sigma_{11},\Sigma_{12},\Sigma_{21})}$ is simply $S_3$. 
\vskip 0.30cm 
This example shows that how the approach works. First we use the 
category of representations of a finite group to find how the tensor products  
of irreducible representations decomposite. Then we find
whether there exist corresponding semi-stable parabolic vector bundles 
which give an equivalent category. As such, the problem becomes a 
combinatoric one, since on ${\Bbb P}^1$,  all vector bundles are direct
sums of ${\Cal O}(n)$ for $n\in {\Bbb Z}$. Further examples  with
all dihedral groups may be easily constructed as well, since their irreducible
representations are at most of rank 2. We leave this to the reader.

Up to this point, naturally, we may ask to which extend this approach to
the Inverse Galois Problem works over number fields $K$. To answer this, 
we come back to moduli spaces of parabolic vector bundles constructed by
Narasimhan and Seshadri. Unlike existing approaches to this problem, 
we have enough rational points ready to use: the anti-canonical line bundles of
these moduli spaces are (roughly speaking) positive. Therefore. for
us, the essential problem becomes
whether the categories generated by $K$-rational parabolic bundles 
correspond to coverings defined over $K$. 
\vfill
\eject
\vskip 0.45cm 
\centerline {\li REFERENCES} 
\vskip 0.25cm 
\item{[DM]} P. Deligne \& J.S. Milne, Tannakian categories, in {\it Hodge Cycles, Motives and 
Shimura Varieties}, LNM {\bf 900}, (1982), 101-228 
\vskip 0.25cm 
\item{[Ka]} E. van Kampen, Almost periodic functions and compact groups, 
Ann. of Math. {\bf 37} (1936), 78-91 
\vskip 0.25cm 
\item{[MS]} V.B.  Mehta \& C.S.  Seshadri, Moduli of vector bundles on curves with 
parabolic structures.  Math. Ann.  {\bf 248}  (1980), no. 3, 205--239. 
\vskip 0.25cm 
\item{[MFK]} D. Mumford, J. Fogarty \& F. Kirwan, {\it Geometric 
Invariant Theory}, Springer-Verlag, 1994 
\vskip 0.25cm 
\item{[Na]} M.  Namba, {\it Branched coverings and algebraic functions}. 
Pitman Research Notes in Mathematics Series {\bf 161}, Longman Scientific 
\& Technical, 1987 
\vskip 0.25cm 
\item{[NS]} M.S.  Narasimhan \& C.S.  Seshadri,  Stable 
and unitary vector bundles on a compact Riemann surface. Ann. of Math. 
{\bf 82} (1965), 540-567 
\vskip 0.25cm 
\item{[Se1]} C.S. Seshadri,  Moduli of $\pi$-vector bundles over an 
algebraic curve.  {\it Questions on Algebraic Varieties}, C.I.M.E., 
III, (1969) 139--260 
\vskip 0.25cm 
\item {[Se2]} C. S. Seshadri, {\it Fibr\'es vectoriels sur les courbes alg\'ebriques}, Asterisque 
{\bf 96}, 1982 
\vskip 0.25cm 
\item{[Ta]} T. Tannaka, {\it Theory of topological groups}, Iwanami, 1949 
(in Japanese) 
\vskip 0.25cm 
\item{[W]} A. Weil, G\'en\'eralisation des fonctions ab\'eliennes, J. 
Math Pures et Appl, {\bf 17}, (1938) 47-87 
\vskip 0.25cm 
\item{[We]} L. Weng, A Program for Gromatric Arithmetic, at math.AG/011240 
\end